\documentclass[english, 12pt,leqno]{amsart}

\usepackage{amsmath}
\usepackage{amsfonts}
\usepackage{amssymb}

\newtheorem{Def}{D\'efinition}

\newtheorem{Cor}[Def]{Corollary}

\newtheorem{Lemma}[Def]{Lemma}
\newtheorem{Prop}[Def]{Proposition}

\newtheorem{Thm}[Def]{Theorem}

\title[Faithful irreducible projective representations]
{Groups with faithful irreducible projective unitary representations}

\author[Bachir Bekka and Pierre de la Harpe]
{Bachir Bekka and Pierre de la Harpe}

\date{September 17, 2010}
\thanks{The authors acknowledge support 
from the Swiss national science foundation}

\address{
B. Bekka: 
UFR Math\'ematiques, 
Universit\'e de Rennes~1, 
Campus Beaulieu,
F--35042 Rennes Cedex.
}
\email{bachir.bekka@univ-rennes1.fr}

\address{
P. de la Harpe:
Section de math\'ematiques, 
Universit\'e de Gen\`eve, 
C.P.~64, 
CH--1211 Gen\`eve 4. 
}
\email{Pierre.delaHarpe@unige.ch}


\begin{document}

\begin{abstract} 
For a countable group $\Gamma$ and a multiplier
$\zeta \in Z^2(\Gamma, \mathbf T)$,
we study the property of $\Gamma$  having a 
unitary projective $\zeta$-repre\-sentation
which is both irreducible and projectively faithful.
Theorem \ref{Thm1} shows that
this property is equivalent to $\Gamma$ being
the quotient  of an appropriate group by its centre.
Theorem~\ref{Thm4} gives a criterion in terms
of the minisocle of $\Gamma$.
Several examples are described
to show the existence of various behaviours.
\end{abstract}

\maketitle

\section{\textbf{
Introduction}}
\label{section1}
% s1

For a Hilbert space $\mathcal H$, we denote by
$\mathcal U (\mathcal H)$ the group of its unitary operators.
We identify
$\mathbf T :=  \left\{ z \in \mathbf C 
\hskip.1cm \big\vert \hskip.1cm
\vert z \vert = 1 \right\}$
with the centre of $\mathcal U (\mathcal H)$, 
namely with the scalar multiples 
of the identity operator $\operatorname{id}_{\mathcal H}$,
we denote by
$\mathcal P \mathcal U (\mathcal H) :=  \mathcal U (\mathcal H) / \mathbf T$ 
the \emph{projective unitary group} of $\mathcal H$, and by
\begin{equation*}
p_{\mathcal H} \,  :  \,
\mathcal U (\mathcal H) \longrightarrow \mathcal P \mathcal U (\mathcal H)
\end{equation*}
the canonical projection.

Let $\Gamma$ be a group.
A projective unitary representation,
or shortly here a \emph{projective representation},
of $\Gamma$ in $\mathcal H$  is a mapping
\begin{equation*}
\pi : \Gamma \,  \longrightarrow \,  \mathcal U (\mathcal H)
\end{equation*}
such that $\pi(e) = \operatorname{id}_{\mathcal H}$
and such that the composition
\begin{equation*}
\underline{\pi} := p_{\mathcal H}\pi \, : \, 
\Gamma \,  \longrightarrow \,  \mathcal P \mathcal U (\mathcal H) 
\end{equation*}
is a homomorphism of groups.
When we find it useful, we write $\mathcal H_{\pi}$ for
the Hilbert space of a projective representation $\pi$.
\par

A projective representation $\pi$ of a group $\Gamma$ is 
\emph{projectively faithful}, or shortly \emph{P-faithful},
if the corresponding homomorphism $\underline{\pi}$ is injective.
The \emph{projective kernel} of  $\pi$ 
is the normal subgroup
\begin{equation}
\label{eq1}
% eq1
\operatorname{Pker}(\pi) \, = \, \operatorname{ker}(\underline{\pi}) \, = \, 
\left\{ x \in \Gamma
\hskip.1cm \big\vert \hskip.1cm
\pi(x) \in \mathbf T \right\} 
\end{equation}
of $\Gamma$. 
In case $\pi$ is a unitary representation,
$\ker (\pi)$ is a subgroup of $\operatorname{Pker}(\pi)$,
sometimes called the the \emph{quasikernel}  of $\pi$,
which can be a proper subgroup,
so that faithfulness of $\pi$ does not imply P-faithfulness.
% even if $\pi$ happens to be a homomorphism.
\par

A projective representation $\pi$  is \emph{irreducible}
if the only closed $\pi(\Gamma)$-invariant subspace of $\mathcal H_{\pi}$ 
are $\{0\}$ and $\mathcal H_{\pi}$.
\par

As a continuation of \cite{BeHa--08}, the present paper results of
our effort to understand
which groups have irreducible P-faithful projective representations.
Our first observation is a version in the present context of
Satz 4.1 of \cite{Pahl--68}. We denote by $Z(\Gamma)$ the centre
of a group $\Gamma$.

\begin{Thm}
\label{Thm1}
% T1
For a group $\Gamma$, the following two properties are equivalent:
\begin{itemize}
\item[(i)]
$\Gamma$ affords an irreducible P-faithful  projective representation;
\item[(ii)]
there exists a group $\Delta$ 
which affords an irreducible faithful unitary representation  
and which is such that 
$\Delta / Z(\Delta) \approx \Gamma$.
\end{itemize}
If, moreover, $\Gamma$ is countable, these properties are also equivalent to:
\begin{itemize}
\item[(iii)]
there exists a {\rm countable} group $\Delta$ 
as in (ii).
\end{itemize}
\end{Thm}

Countable groups which have irreducible faithful unitary representations 
have been characterised in \cite{BeHa--08}, 
building up on results of \cite{Gasc--54} for finite groups.
\par

A group $\Gamma$ is \emph{capable} if there exists a group $\Delta$
with $\Gamma \approx \Delta / Z(\Delta)$,
and \emph{incapable} otherwise. 
The notion appears in \cite{Baer--38},
which contains a criterion of capability 
for abelian groups which are direct sums of cyclic groups 
(for this, see also \cite{BeFS--79}),
% (see also Corollary 4.10, Page 219 of \cite{BeTa--82}),
and the terminology ``capable'' is that of \cite{HaSe--64}.
Conditions for capability (several of them being either necessary or sufficient)
are given in Chapter IV of \cite{BeTa--82}.
\par

The \emph{epicentre} of a group $\Gamma$ is the largest central subgroup $A$
such that the quotient projection $\Gamma \longrightarrow \Gamma/A$
induces in homology an injective homomorphism 
$H_2(\Gamma,\mathbf Z) \longrightarrow H_2(\Gamma/A, \mathbf Z)$,
where $\mathbf Z$ is viewed as a trivial module. 
This group was introduced in \cite{BeFS--79} 
and \cite{BeTa--82}, with a formally different definition;  
the terminology is from \cite{Elli--98},
and the characterisation given above appears in Theorem 4.2 of \cite{BeFS--79}.

\begin{Prop}[Beyl-Felgner-Schmid-Ellis]
\label{PropBFSE}
% Prop 2
Let $\Gamma$ be a group and let $Z^*(\Gamma)$ denote its epicentre.
\par

(i) $\Gamma$ is capable if and only if $Z^*(\Gamma) = \{e\}$.
\par

(ii) $\Gamma / Z^*(\Gamma)$ is capable in all cases.
% $\Gamma / Z^*(\Gamma)$ est capable \cite{BeTa--82}, page 209.

(iii) A perfect group with non trivial centre is incapable. 

\end{Prop}

\begin{Cor}
\label{Cor3}
A perfect group with non-trivial centre has no P-faithful projective representation.
\end{Cor}

A \emph{multiplier} on $\Gamma$ is a mapping 
$\zeta : \Gamma \times \Gamma \longrightarrow \mathbf T$
such that
\begin{equation}
\label{eq2}
% eq2
\zeta(e,x) \, = \, \zeta(x,e) \, = \, 1 
\hskip.3cm \text{and} \hskip.3cm
\zeta(x,y)\zeta(xy,z) \, = \, \zeta(x,yz)\zeta(y,z)
\end{equation}
for all $x,y,z \in \Gamma$.
We denote by $Z^2(\Gamma, \mathbf T)$ the set of all these,
which is an abelian group for the pointwise product.
A projective representation $\pi$ of $\Gamma$ in $\mathcal H$
determines a unique multiplier $\zeta_{\pi} = \zeta$ such that
\begin{equation}
\label{eq3}
% eq3
\pi(x) \pi(y) \, = \, \zeta(x,y) \pi(xy)
\end{equation}
for all $x,y \in \Gamma$; we say then that $\pi$ is a 
\emph{$\zeta$-representation} of $\Gamma$.
Conversely, any $\zeta \in Z^2(\Gamma, \mathbf T)$ occurs in such a way;
indeed, $\zeta$ is the multiplier determined by the
\emph{twisted left regular $\zeta$-representation} 
$\lambda_{\zeta}$ of $\Gamma$, 
defined on the Hilbert space $\ell^2(\Gamma)$ by
\begin{equation}
\label{eq4}
% eq4
\left( \lambda_{\zeta}(x) \varphi \right) (y) \, = \,
\zeta(x, x^{-1}y) \varphi(x^{-1}y) .
\end{equation}
A good reference for these regular $\zeta$-representations is \cite{Klep--62}.
In the special case $\zeta = 1$, a $\zeta$-representation of $\Gamma$
is just a \emph{unitary representation} of $\Gamma$; but we repeat that 

\medskip

\noindent \textbf{First standing assumption.} 
In this paper, \emph{by ``representation'', 
we always mean ``unitary representation''.}

\medskip

For a projective representation of $\Gamma$ in $\mathbf C$,
namely for a mapping $\nu : \Gamma \longrightarrow \mathbf T$ with $\nu(e) = 1$,
let $\zeta_{\nu} \in Z^2(\Gamma, \mathbf T)$
denote the corresponding multiplier, namely the mapping defined by
\begin{equation}
\label{eq5}
% eq5
\zeta_{\nu}(x,y) \, = \,  \nu(x) \nu(y) \nu(xy)^{-1} .
\end{equation}
We denote by $B^2(\Gamma, \mathbf T)$ the set of 
all multipliers of the form $\zeta_{\nu}$,
which is a subgroup of $Z^2(\Gamma, \mathbf T)$,
and by 
$H^2(\Gamma, \mathbf T) := 
Z^2(\Gamma, \mathbf T) / B^2(\Gamma, \mathbf T)$
the quotient group;
as usual, $\zeta,\zeta' \in Z^2(\Gamma, \mathbf T)$ are \emph{cohomologous}
if they have the same image in $H^2(\Gamma, \mathbf T)$.
\par

Given a $\zeta$-representation $\pi$ of $\Gamma$ in $\mathcal H$,
there is a standard bijection between:
\begin{itemize}
\item[$\circ$] 
the set of projective representations 
$\pi' : \Gamma \longrightarrow \mathcal U (\mathcal H)$
such that $p_{\mathcal H} \pi = p_{\mathcal H} \pi'$, on the one hand,
\item[$\circ$] 
and the set of multipliers  cohomologous to $\zeta$, on the other hand.
\end{itemize}
In other terms, a group homomorphism $\underline \pi$ of $\Gamma$
in $\mathcal P \mathcal U (\mathcal H)$
determines a class\footnote{This is of course the class associated
to the extension of $\Gamma$ by $\mathbf T$ pulled back by $\underline \pi$
of the extension
$\{1\} \longrightarrow \mathbf T \longrightarrow \mathcal U (\mathcal H)
\longrightarrow \mathcal P \mathcal U (\mathcal H) \longrightarrow \{1\}$.
See for example Section 6.6 in \cite{Weib--94}.} 
in $H^2(\Gamma, \mathbf T)$,
and the set of projective representations covering $\underline \pi$
is in bijection with the representatives of this class in $Z^2(\Gamma, \mathbf T)$.
Observe that $\pi$ and $\pi'$  above are together irreducible or not,
and together P-faithful or not. 

\par

For much more on projective representations and multipliers,
in the setting of separable locally compact groups,
see \cite{Mack--58};
for a very short but informative exposition on earlier work,
starting with that of Schur, see \cite{Kall--84}.
Note that other authors (such as Kleppner) use ``multiplier representation''
and ``projective representation'' 
when Mackey  uses  ``projective representation'' 
and ``homomorphism in $\mathcal P \mathcal U (\mathcal H)$'',
respectively.

\medskip

\noindent \textbf{Definition.} 
\emph{
Given a group $\Gamma$ and a multiplier $\zeta \in Z^2(\Gamma, \mathbf T)$, 
the group $\Gamma$ is}
irreducibly $\zeta$-represented
\emph{if it has an irreducible P-faithful  $\zeta$-representation.}
\par

This depends only on the class 
$\underline{\zeta} \in H^2(\Gamma, \mathbf T)$ 
of $\zeta$.

\medskip

For a group $\Gamma$, recall that 
a \emph{foot} of $\Gamma$ is a minimal normal subgroup,
that the \emph{minisocle} is the subgroup $MS(\Gamma)$ of $\Gamma$
generated by the union of all finite feet of $\Gamma$, 
and that $MA(\Gamma)$ is the subgroup of $MS(\Gamma)$
generated by the union of all finite abelian feet of $\Gamma$.
It is obvious that $MS(\Gamma)$ and $MA(\Gamma)$ 
are characteristic subgroups of $\Gamma$;
it is easy to show that $MA(\Gamma)$ is abelian
and is a direct factor of $MS(\Gamma)$.
For all  this, we refer to  Proposition 1 in  \cite{BeHa--08}.

Let $N$ be a normal subgroup of $\Gamma$ 
and $\sigma$  a $\zeta$-representation of $N$, 
for some $\zeta \in Z^2(N, \mathbf T)$.
If $\zeta = 1$ (the case of ordinary representations), 
define the \emph{$\Gamma$-kernel of $\sigma$} by
\begin{equation*}
\ker_{\Gamma}(\sigma) \, = \, 
\ker \Big( \bigoplus_{\gamma \in \Gamma} \sigma^{\gamma} \Big)
\end{equation*}
where $\sigma^{\gamma}(x) := \sigma(\gamma x \gamma^{-1})$;
say, as in \cite{BeHa--08}, that $\sigma$ is \emph{$\Gamma$-faithful}
if this $\Gamma$-kernel is reduced to $\{e\}$;
when $\zeta$ is the restriction to $N$ of a multiplier 
(usually denoted by $\zeta$ again)
in $Z^2(\Gamma, \mathbf T)$,
there is an analogous notion for the general case ($\zeta \ne 1$), 
called \emph{$\Gamma$-P-faithfulness},
used in Theorem \ref{Thm4}, 
but defined only in Section \ref{section3} below.
Before the next result, 
we find it convenient to define one more property.

\medskip

\noindent \textbf{Definition.} 
\emph{
A group $\Gamma$ has
} 
Property (Fab)
\emph{
if any normal subgroup of $\Gamma$ generated by one conjugacy class
has a finite abelianisation.}

\medskip

Examples of groups which enjoy Property (Fab) include
finite groups, $\operatorname{SL}_n(\mathbf Z)$ for $n \ge 3$,
and more generally 
lattices in a finite product $\prod_{\alpha \in A} G_{\alpha}$
of simple groups $G_{\alpha}$ 
over (possibly different) local fields $k_{\alpha}$
when $\sum_{\alpha \in A} k_{\alpha}-{\rm rank} (G_{\alpha}) \geq 2$ 
(see \cite{Marg--91}, IV.4.10, and Example VI below).
They also include abelian locally finite groups,
and more generally torsion groups which are FC, namely which are such
that all their conjugacy classes are finite;
in particular, they include groups of the form $MS(\Gamma)$ and $MA(\Gamma)$.

\begin{Thm}
\label{Thm4}
Let $\Gamma$ be a countable group and let $\zeta \in Z^2(\Gamma, \mathbf T)$.
Consider the following conditions:
\begin{itemize}
\item[(i)]
$\Gamma$ is irreducibly $\zeta$-represented;
\item[(ii)]
$MS(\Gamma)$ has a $\Gamma$-P-faithful irreducible $\zeta$-representation;
\item[(iii)]
$MA(\Gamma)$ has a $\Gamma$-P-faithful irreducible $\zeta$-representation.
\end{itemize}
\par\noindent
Then $(i) \Longrightarrow (ii) \Longleftrightarrow (iii)$.
\par
If, moreover, $\Gamma$ has Property (Fab), then $(ii) \Longrightarrow (i)$,
so that $(i)$, $(ii)$, and $(iii)$ are equivalent.
\end{Thm}

The hypothesis ``$\Gamma$ countable''  is essential 
because our arguments use measure theory and direct integrals; 
in fact, Theorem \ref{Thm4} fails in general for uncountable groups 
(see Example (VII) in \cite{BeHa--08}, Page 863). 
About the converse of $(ii) \Longrightarrow (i)$,
see Example I below.
\par

Recall that a group $\Gamma$ has \emph{infinite conjugacy classes}, or is \emph{icc},
if $\Gamma \ne \{e\}$ and if any conjugacy class in $\Gamma \smallsetminus \{e\}$
is infinite.
For example, a lattice in a 
centreless connected semisimple Lie group without compact factors 
is icc, as a consequence of Borel Density Theorem
(see Example VI).

\begin{Cor}
\label{Cor5}
Let $\Gamma$ be a countable group which has Property (Fab)
and which fulfills at least one of the three following conditions:
\begin{itemize}
\item[(i)]
$\Gamma$ is torsion free;
\item[(ii)]
$\Gamma$ is icc;
\item[(iii)]
$\Gamma$ has a faithful primitive action on an infinite set.
\end{itemize}
Then, for any $\zeta \in Z^2(\Gamma, \mathbf T)$,
the group $\Gamma$ is irreducibly $\zeta$-represented.
\end{Cor}

Indeed, any of  Conditions (i) to (iii) implies that
$MS(\Gamma) = \{e\}$. 
Recall that, if $\Gamma$ fulfills (iii) on an infinite set $X$, 
any normal subgroup $N \ne \{e\}$ acts 
transitively on $X$, and therefore  is infinite
(see \cite{GeGl--08}).

\medskip

A group can be either irreducibly represented or not,
and also either irreducibly  $\zeta$-represented or not (for some $\zeta$).
These dichotomies separate groups in four classes,
each one illustrated in Section \ref{sectionExamples}
by one of Examples I to IV below. 
Examples V and VI illustrate the same class as Example~I. 
\par

Section \ref{section3} contains standard material on multipliers,
and the definition of $\Gamma$-P-faithfulness;
mind the ``second standing assumption'' on the normalisation of multipliers, 
which applies to all other sections.
In Section \ref{section4}, we review central extensions
and prove Theorem \ref{Thm1}.
Sections \ref{section5} and \ref{section6} contain the proof
of Theorem \ref{Thm4}, respectively the part which involves
our ``Property (Fab)'' and the part which does not.
\par

Section \ref{section7} contains material on (in)capability,
and the proof of Proposition \ref{PropBFSE}.
Section \ref{section8} describes a construction of 
irreducible P-faithful projective representation of  a class of abelian groups,
and expands on Example II.
The last section is a digression to point out
a fact from homological algebra
which in our opinion is not quoted often enough in the literature
on projective representations.

\section{\textbf{
Examples}}
\label{sectionExamples}
% section2

\noindent \textbf{Example I.}
The implication $(ii) \Longrightarrow (i)$ of Theorem \ref{Thm4}
does not hold for a free abelian group $\mathbf Z^n$ ($n \ge 1$)
and the unit multiplier\footnote{
Recall that $H^2(\mathbf Z, \mathbf T)= \{0\}$,
because $\mathbf Z$ is free, 
% Weibel, corollary 6.2.7
so that
\begin{equation*}
Z^2(\mathbf Z, \mathbf T) = B^2(\mathbf Z, \mathbf T) =
\operatorname{Mapp}(\mathbf Z, \mathbf T) / 
\operatorname{Hom}(\mathbf Z, \mathbf T) .
\end{equation*}
Also $H^2(\mathbf Z^n, \mathbf T)= \mathbf Z^{n(n-1)/2}$ for all $n \ge 1$.
}
$\zeta = 1 \in Z^2(\mathbf Z^n, \mathbf T)$.
Indeed, on the one hand, Condition (ii) of Theorem \ref{Thm4}
is satisfied since $MS(\mathbf Z^n) = \{0\}$.
On the other hand, since $\mathbf Z^n$ is abelian, 
any irreducible $\zeta$-representation
(that is any ordinary irreducible representation)
is one-dimen\-sional, so that its projective kernel is the whole of $\mathbf Z^n$,
and therefore $\mathbf Z^n$ is \emph{not} irreducibly $\zeta$-represented.
Moreover, $\mathbf Z^n$ being for any $n \ge 1$ a (dense) subgroup of $\mathbf T$,
it has an irreducible faithful representation of dimension one.

\medskip

\noindent \textbf{Example II.}
There are groups which do not afford any irreducible faithful representation
but which do have projective representations which are irreducible and P-faithful.
\par

The Vierergruppe 
$\mathbf V = \mathbf Z / 2\mathbf Z \times \mathbf Z / 2\mathbf Z$,
being finite abelian non-cyclic,
does not have any irreducible faithful representation.
If $\zeta \in Z^2(\mathbf V, \mathbf T)$ is  a cocycle representing
the non-trivial cohomology class in 
$H^2(\mathbf V, \mathbf T) \approx \mathbf Z / 2\mathbf Z$, 
then $\mathbf V$ has a $\zeta$-representation of degree $2$ 
which is both irreducible and P-faithful,
essentially given by the Pauli matrices 
(see Section IV.3 in \cite{Simo--96}).

Part of this carries over to any non-trivial finite abelian group of the form
$L \times L$. More on this in Section \ref{section8}.

\medskip

\noindent \textbf{Example III.}
Let us first recall a few basic general facts
about  irreducible projective representations of a finite group $\Gamma$.
The cohomology group 
$H^2(\Gamma,\mathbf T)$ is isomorphic to the homology group
$H_2(\Gamma, \mathbf Z)$, and is finite.
Choose a multiplier $\zeta \in Z^2(\Gamma, \mathbf T)$,
say normalised (see Section \ref{section3} below).
An element $x \in \Gamma$ is \emph{$\zeta$-regular} if
$\zeta(x,y) = \zeta(y,x)$ whenever $y \in \Gamma$ commutes with $x$;
it can be checked that a conjugate of a regular element is again regular.
Let $h(\zeta)$ denote the number of  conjugacy classes
of $\zeta$-regular elements in $\Gamma$.
Then it is known that $\Gamma$ has exactly $h(\zeta)$
irreducible $\zeta$-representations, up to unitary equivalence,
say of degrees  $d_1, \hdots, d_{h(\zeta)}$;
moreover each $d_j$ divides the order of $\Gamma$, 
and $\sum_{j=1}^{h(\zeta)} d_j^2 = \vert \Gamma \vert$.
See Chapter 6 in \cite{BeZh--98}, in particular Corollary 10 and Theorem 13 Page 149.
\par

Clearly $h(\zeta) \le h(1)$ for all $\zeta \in Z^2(\Gamma, \mathbf T)$.
It follows from Lemma \ref{lemma11} below that, if $\underline{\zeta} \ne 1$,
then $d_j \ge 2$ for all $j \in \{1, \hdots, h(\zeta)\}$.
\par

Now, for the gist of this Example III, 
assume  that $\Gamma$ is  a nonabelian finite simple group.
Then, except for the unit character, 
any representation of $\Gamma$
is faithful and any projective representation of $\Gamma$ is P-faithful.

\medskip

\noindent \textbf{Example IV.}
Let $\Gamma$ be a perfect group.
Its \emph{universal central extension} $\widetilde \Gamma$
is a perfect group with centre the\footnote{Recall that, 
for a perfect group $\Gamma$
and a trivial $\Gamma$-module $A$, we have
$H^2(\Gamma,A) \, \approx \, 
\operatorname{Hom}(H_2(\Gamma,\mathbf Z),A)$
as a consequence of the
universal coefficient theorem for cohomology;
in particular, $H^2(\Gamma, \mathbf T)$
is the Pontryagin dual of the Schur multiplier.
In case $H_2(\Gamma, \mathbf Z)$ is moreover finite, 
e.g. if $\Gamma$ is perfect and finite,
$H^2(\Gamma, \mathbf T)$ is isomorphic to the Schur multiplier
(non-canonically).}
\emph{Schur multiplier}
$H_2(\Gamma,\mathbf Z)$
and central quotient $\Gamma$ \cite{Kerv--70}.
If this Schur multiplier is not $\{0\}$, 
$\widetilde \Gamma$ is incapable,
and therefore does not have any irreducible P-faithful projective representation
(Corollary \ref{Cor3}).
% by Proposition \ref{PropBFSE}.
If $\Gamma$ is as in (i) or (ii) below,
$\widetilde \Gamma$ is moreover
not irreducibly represented (by \cite{Gasc--54} and \cite{BeHa--08}):
\begin{itemize}
\item[(i)]
$\Gamma$ is a finite simple group 
with $H_2(\Gamma, \mathbf Z)$ not cyclic\footnote{If
$H_2(\Gamma, \mathbf Z)$ is cyclic not $\{0\}$, see Example V.}.
The complete list of such groups is given
in Theorem 4.236, Page 301 of \cite{Gore--82},
and includes the finite simple group $\operatorname{PSL}_3(\mathbf F_4)$,
also denoted by $A_2(4)$, one of the two finite simple groups of order $20 160$.
\item[(ii)]
$\Gamma$ is one of 
the Steinberg groups $\operatorname{St}_3(\mathbf Z)$ 
and $\operatorname{St}_4(\mathbf Z)$,
which are  the universal central extensions of 
$\operatorname{SL}_3(\mathbf Z)$ 
and $\operatorname{SL}_4(\mathbf Z)$, respectively.
Indeed, 
van der Kallen\footnote{Thanks to Andrei Rapinchuk
for this reference.}
% et pour \cite{Rose--94},  voir la partie (3) de l'exercice 4.3.20 
\cite{Kall--74} has shown that
\begin{equation*}
H_2(\operatorname{SL}_3(\mathbf Z), \mathbf Z) \, \approx \, 
H_2(\operatorname{SL}_4(\mathbf Z), \mathbf Z) \, \approx \, 
\mathbf Z / 2\mathbf Z \times \mathbf Z / 2\mathbf Z .
\end{equation*}
Thus, these groups are not irreducibly represented by \cite{BeHa--08}.
(For $n \ge 5$, it is known that  
$H_2(\operatorname{SL}_n(\mathbf Z), \mathbf Z) \approx \mathbf Z / 2\mathbf Z$;
see \cite{Miln--71}, Page 48.
And $H_2(\operatorname{SL}_2(\mathbf Z), \mathbf Z) = \{0\}$,
see the comments after Proposition \ref{PropIncapables}.)
\end{itemize}

\medskip

\noindent \textbf{Example V.}
This example and the next one will show,
besides  the $\mathbf Z^n$~'s of Example~I, groups
which are irreducibly represented,
but which do not have any irreducible P-faithful representation.

\par

Any finite perfect  group $\Gamma$ with 
centre $Z(\Gamma)$ cyclic and not $\{0\}$ has these properties, 
by Gasch\"utz theorem and by Corollary \ref{Cor3}. 
This is for example the case of the quasi simple group 
$\operatorname{SL}_n(\mathbf F_q)$
whenever the finite field $\mathbf F_q$ has non-trivial $n$th roots of unity,
so that 
$Z(\operatorname{SL}_n(\mathbf F_q)) = 
\{\lambda \in \mathbf F_q \mid \lambda^n = 1 \}$
is cyclic and not $\{e\}$.
(As usual, $\operatorname{SL}_2(\mathbf F_2)$ 
and $\operatorname{SL}_2(\mathbf F_3)$ are ruled out.)
\par

The groups $\operatorname{SL}_{2n}(\mathbf Z)$, for $2n \ge 4$,
are perfect with centre cyclic of order $2$, and therefore incapable,
so that Corollary \ref{Cor3} applies;
the group $\operatorname{SL}_2(\mathbf Z)$, which is not perfect, 
is also incapable (see Section \ref{section7}).
On the other hand, since for all $n \ge 1$
the minisocle of $\operatorname{SL}_{2n}(\mathbf Z)$
coincides with its centre, of order $2$, 
these groups do have representations which are irreducible and faithful.
These considerations hold also for the symplectic groups
$\operatorname{Sp}_{2n}(\mathbf Z)$, $2n \ge 6$,
which are perfect  \cite{Rein--95}.

\medskip 

\noindent \textbf{Example VI.a.}
Let $B$ be a finite set. For $\beta \in B$,
let $\mathbf k_{\beta}$ be a local field 
and $\mathbf G_{\beta}$ be a nontrivial connected semi-simple group
defined over $\mathbf k_{\beta}$,
without $\mathbf k_{\beta}$-anisotropic factor.
Set $G = \prod_{\beta \in B} \mathbf G_{\beta}(\mathbf k_{\beta})$,
with its locally compact topology which makes it
a $\sigma$-compact, metrisable, compactly generated group.
Let $\Gamma$ be an irreducible lattice in $G$.
\par

If $N$ is a finite normal subgroup of $\Gamma$, 
we claim that $N$ is central in $\Gamma$.
If there are several factors ($\vert B \vert \ge 2$), 
the claim is a consequence of the fact that the projection of the lattice
in each factor is dense, by irreducibility. 
If $\vert B \vert = 1$, consider $x \in N$.
The centraliser $Z_{\Gamma}(x)$ of $x$ in $\Gamma$
is also a lattice in $G$ because it is of finite index in $\Gamma$.
By the Borel-Wang density theorem (Corollary 4.4 of Chapter II in \cite{Marg--91}),
$Z_{\Gamma}(x)$ is Zariski-dense in $G$, 
so that $x$ commutes with every element of $G$,
and this proves the claim.
\par

If follows that, if moreover the centre of $G$ if finite cyclic, 
then $MS(\Gamma) = MA(\Gamma)$ is also a finite cyclic group,
so that $\Gamma$ is irreducibly represented by \cite{BeHa--08}. 

\medskip 

\noindent \textbf{Example VI.b.}
To continue this same example, let us particularise the situation to the case
of a non-compact semi-simple real Lie group $G$,
which is connected, not simply connected, and with a non-trivial centre.
Let $\Gamma$ be a lattice in
$G$ with a non-trivial  centre $Z(\Gamma)$.
Denote by $\widetilde \Gamma$ the inverse image of  $\Gamma$ in the universal
cover of $G$ and by $p:\widetilde \Gamma\to \Gamma$ the canonical projection.
Observe that  $Z(\widetilde \Gamma)=p^{-1}(Z(\Gamma))$. 
% Moreover,  the kernel $\ker p$ is contained in  
% $Z(\widetilde \Gamma)$ as a proper subgroup, because
% $G$ is not simply connected and $Z(\Gamma)\neq \{e\}$.
Choose a set-theoretical section $s: \Gamma\to \widetilde\Gamma$ 
for $p$ with $s(e)=e$ 
and a character $\chi\in {\rm Hom}(  Z(\widetilde \Gamma), \mathbf T)$.
Define a mapping
\begin{equation*}
\zeta : \Gamma \times \Gamma \longrightarrow \mathbf T, 
\qquad 
\zeta(x,y) \, = \, \chi \left(s(x)s(y) s(xy)^{-1}\right).
\end{equation*}
It is straightforward to check that $\zeta$ is a multiplier, namely that
$\zeta\in Z^2(\Gamma, \mathbf T)$.

[Classes of multipliers of this kind are not  arbitrary. They correspond precisely to
those classes in  $H^2(\Gamma, \mathbf T)$ which are restrictions
of classes in the appropriately defined group $H^2(G, \mathbf T)$.
The latter group is known to be isomorphic to ${\rm Hom}(\pi_1(G), \mathbf T)$;
see Proposition 3.4 in \cite{Moor--64} and \cite{BaMi--00}.]

 It is obvious that, if $\chi$ extends to a unitary character $\widetilde\chi$
of $\widetilde \Gamma$,  then $\zeta$ belongs to $B^2(\Gamma, \mathbf T)$.
Indeed, in this case $\zeta = \zeta_\nu$ for $\nu : \Gamma \longrightarrow \mathbf T$
defined by $\nu(x)= \widetilde{\chi}(s(x))$.
Conversely, assume that  $\zeta = \zeta_\nu$ for  some
$\nu : \Gamma \longrightarrow \mathbf T$ with $\nu(e)= 1$.
As in the proof of Theorem 1.1 in \cite{BaMi--00},
define $v : \widetilde{\Gamma} \longrightarrow  {\mathbf T}$ by
$v(x) = \chi (x^{-1} s(p(x)))^{-1}$  for $x \in \widetilde\Gamma$
(observe that  $ x^{-1}s(p(x)) \in Z(\widetilde\Gamma)$, so that $v(x)$ is well-defined).
One  checks that
\begin{equation*}
\zeta ( p(x), p(y)) \, = \,  v(x) v(y) v(xy))^{-1} \qquad \forall x,y\in \widetilde \Gamma.
\end{equation*}  
It follows that the function 
$\widetilde{\chi} : \widetilde{\Gamma} \longrightarrow {\mathbf T}$,
defined  by $ \widetilde{\chi} (x) =  \nu (p(x)) v(x)^{-1}$, is a character
of $\widetilde \Gamma$ which extends $\chi$.

As a consequence,  if the intersection of $Z(\widetilde \Gamma)$ with  the
commutator subgroup $[\widetilde \Gamma, \widetilde \Gamma]$ is
not reduced to $\{e\}$, we can find $\chi$ such that the corresponding multiplier
$\zeta$ does not belong to  $B^2(\Gamma, \mathbf T)$.
We provide in VI.c below an example for which this does occur.

We claim that $\Gamma$ has no P-faithful irreducible $\zeta$-representation.
By Theorem 4, it suffices to show that  $MS(\Gamma)$ has no $\Gamma$-P-faithful
irreducible   $\zeta$-representation.

The group $MS(\Gamma)$ coincides with  $Z(\Gamma)$ as we have shown 
in Part VI.a of the present example. 
Observe  that $s(x)\in Z(\widetilde \Gamma)$ for every $x\in Z(\Gamma)$. 
We have therefore $\zeta(x,y)=\zeta(y,x)$ for all $x\in  Z(\Gamma)$ and $y \in \Gamma$.
In particular, it follows that the restriction of $\zeta$ to $Z(\Gamma)$ is trivial
(see Lemma 7.2 in \cite{Klep--65}). 
Upon changing $\zeta$ inside its cohomology class,
we can assume that $\zeta(x,y)=1$ for all $x,y\in Z(\Gamma)$.

Let $\sigma$ be an irreducible $\zeta$-representation of $Z(\Gamma)$.
Since the restriction of $\zeta$ to  $Z(\Gamma)$  is trivial, 
we have  P$\ker \sigma = Z(\Gamma)$.
 From the fact that $\zeta(x,y) = \zeta(y,x)$ for all $x \in Z(\Gamma)$ and $y \in \Gamma$,
 it follows that  P$\ker_\Gamma \sigma=Z(\Gamma)$; see Remark (B) after Proposition 8.
Hence, $\sigma$ is not  $\Gamma$-P-faithful
since $Z(\Gamma)$ is non-trivial by assumption.

\medskip 

\noindent \textbf{Example VI.c.}
Let $\Delta$ be the fundamental group  of  a closed surface of genus $2$, 
viewed as a subgroup of $PSL_2(\mathbf R)$.
Let  $\Gamma$ be the inverse image  of $\Delta$
in $SL_2(\mathbf R)$; observe that $Z(\Gamma)$ is the two-element group.
The group $\widetilde \Gamma$,
the discrete subgroup of the universal cover of $SL_2(\mathbf R)$ defined in VI.b,
has a presentation with (see IV.48 in \cite{Harp--00})
\begin{itemize}
\item[]
generators: $a_1, a_2, b_1, b_2, c$
\item[]
and relations: $c$ is central, and $[a_1,b_1][a_2,b_2]=c^2$.
\end{itemize}
In particular, the intersection of $Z(\widetilde \Gamma)$ 
with  $[\widetilde \Gamma, \widetilde \Gamma]$ is non trivial .

\newpage

\section{\textbf{
$\Gamma$-P-faithfulness for projective representations 
\\
of normal subgroups of $\Gamma$}}
\label{section3}
% s3

Let $N$ be a normal subgroup of a group $\Gamma$ 
and let $\zeta \in Z^2(N, \mathbf T)$.
Let $\sigma$ be a $\zeta$-representation of $N$.
\par

For $\gamma \in \Gamma$, the mapping 
\begin{equation*}
{N \ni x \,  \longmapsto \, 
\sigma(\gamma x\gamma^{-1}) \in \mathcal U (\mathcal H_{\sigma})}
\end{equation*}
is in general not a $\zeta$-representation of $N$,
but is a $\zeta^\gamma$-representation of $N$,
where $\zeta^\gamma \in Z^2(N, \mathbf T)$ is defined by
\begin{equation}
\label{eq6}
% eq6
\zeta^\gamma(x,y) \, = \, \zeta(\gamma x\gamma^{-1}, \gamma y\gamma^{-1})
\end{equation}
for all $x,y \in N$.
\par

Suppose moreover that $\zeta$ is the restriction to $N$ 
of some multiplier on $\Gamma$.
Then the multiplier $\zeta^\gamma$ is cohomologous to $\zeta$; more precisely:

\begin{Lemma}[Mackey]
\label{lemma6}
% lemma 6
Let $\zeta \in Z^2(\Gamma,\mathbf T)$, 
let  $\sigma$ be a $\zeta$-representation of $N$,
and let $\gamma \in \Gamma$.

(i)
Define a mapping $\nu_\gamma : N \longrightarrow \mathbf T$ by
$\nu_\gamma(x) =  
\frac{ \zeta(\gamma,x)  \zeta(\gamma x,\gamma^{-1}) }{ \zeta(\gamma^{-1},\gamma) }$.
Then
\begin{equation}
\label{eq7}
% eq7
\frac{ \zeta^\gamma(x,y) }{ \zeta(x,y) } \, = \, 
\frac{ \nu_\gamma(xy) }{ \nu_\gamma(x) \nu_\gamma(y) }
\end{equation}
for all $x,y \in N$.
In particular, $\underline{\zeta^{\gamma}} = \underline{\zeta}
\in H^2(\Gamma, \mathbf T)$.

(ii)
Define a mapping 
$\sigma^\gamma : N \longrightarrow \mathcal U (\mathcal H_{\sigma})$  by
\begin{equation}
\label{eq8}
% eq8
\sigma^\gamma(x) \, = \, 
\zeta(\gamma,x) \zeta(\gamma x,\gamma^{-1}) \, \sigma(\gamma x\gamma^{-1}) .
\end{equation}
Then
\begin{equation}
\label{eq9}
% eq9
\sigma^\gamma(x) \sigma^\gamma(y) \, =  \, 
\zeta(\gamma^{-1},\gamma) \zeta(x,y) \,  \sigma^\gamma(xy)
\end{equation}
for all $x,y \in N$.
\end{Lemma}

\emph{Proof.}
For (i), we refer to Lemma 4.2 in \cite{Mack--58},
of which the proof uses   (\ref{eq2}) from Section \ref{section1}.
For (ii), we have
\begin{equation*}
\aligned
& \sigma^\gamma(x) \sigma^\gamma(y) 
= 
\\
& \hskip.5cm
\zeta(\gamma,x) \zeta(\gamma x,\gamma^{-1}) 
\zeta(\gamma,y) \zeta(\gamma y,\gamma^{-1}) 
\zeta(\gamma x\gamma^{-1},\gamma y\gamma^{-1})  \,  
\sigma(\gamma xy\gamma^{-1}) 
=
\\
& \hskip.5cm
\frac{ 
\zeta(\gamma,x) \zeta(\gamma x,\gamma^{-1}) 
\zeta(\gamma,y) \zeta(\gamma y,\gamma^{-1}) 
}{
\zeta(\gamma,xy) \zeta(\gamma xy,\gamma^{-1}) 
}
\, \frac{ \zeta^\gamma(x,y) }{ \zeta(x,y) } 
\, \zeta(x,y) \, \sigma^\gamma(xy)
=
\\
& \hskip.5cm
\zeta(\gamma^{-1},\gamma) 
\frac{
\nu_\gamma(x) \nu_\gamma(y)
}{
\nu_\gamma(xy)
}
\, \frac{ \zeta^\gamma(x,y) }{ \zeta(x,y) } \, \zeta(x,y) \, \sigma^\gamma(xy)
=
\\
& \hskip.5cm
\zeta(\gamma^{-1},\gamma) \zeta(x,y) \, \sigma^\gamma(xy) ,
\endaligned
\end{equation*}
where we have used (i) in the last equality.
\hfill $\square$

\medskip

Equation (\ref{eq9}) makes it convenient to restrict the discussion
to normalised multipliers.

\medskip

\noindent \textbf{Definition.} 
\emph{A multiplier $\zeta$ on a group $\Gamma$ is} 
normalised
\emph{if $\zeta(x,x^{-1}) = 1$ for all $x \in \Gamma$.}
\emph{A projective representation $\pi$ of a group $\Gamma$ is} 
normalised
\emph{if $\pi(x^{-1}) = \pi(x)^{-1}$ for all $x \in \Gamma$.}

(Some authors, see e.g. Page 142 of \cite{BeZh--98},
use ``normalised'' for multipliers in a different meaning.)

\begin{Lemma}
\label{lemma7}
% lemmma 7
(i) Any multiplier $\zeta'$ on a group 
is cohomologous to a normalised mutliplier $\zeta$.

(ii) If $\zeta$ is a normalised multiplier on a group $\Gamma$, then
\begin{equation}
\label{eq10}
% eq10
\zeta(y^{-1},x^{-1}) \, = \, \zeta(x,y)
\end{equation}
for all $x,y \in \Gamma$.
\end{Lemma}

\emph{Proof.}
(i) 
Let $\pi'$ be an arbitrary $\zeta'$-representation of $\Gamma$
on a Hilbert space $\mathcal H$.
Define $J = \{ \gamma \in \Gamma \hskip.1cm \vert \hskip.1cm \gamma^2 = e \}$
and choose a partition $\Gamma = J \sqcup K \sqcup L$
such that $\ell \in L$ if and only if $\ell^{-1} \in K$.
For each $\gamma \in J$, choose $z_{\gamma} \in \mathbf T$
such that $(z_{\gamma} \pi'(\gamma))^2 = \operatorname{id}_{\mathcal H}$.
Define a projective representation $\pi$ of $\Gamma$ on $\mathcal H$
by $\pi(\gamma) = z_{\gamma} \pi'(\gamma)$ if $\gamma \in J$,
by $\pi(\gamma) = \pi'(\gamma)$ if $\gamma \in K$,
and by $\pi(\gamma) = \pi'(\gamma^{-1})^{-1}$ if $\gamma \in L$.
Then $\pi$ is a normalised projective representation
of which the multiplier $\zeta$ is normalised and cohomologous to $\zeta'$.
\par

(ii) 
This is a  consequence of the identity
\begin{equation*}
\pi(x^{-1}) \pi(y^{-1}) \, = \,  \frac{1}{\zeta(x,y)} \pi(y^{-1}x^{-1}) ,
\end{equation*}
which is a way of writing Equation (\ref{eq3}) when $\pi$ is normalised.
\hfill $\square$

\begin{Prop}
\label{Prop8}
% Proposition 8
Let $\Gamma$ be a group, 
$\zeta \in Z^2(\Gamma, \mathbf T)$  a normalised multiplier,
and $N$ a normal subgroup of $\Gamma$.
% we denote restriction of $\zeta$ to $N$ by the same letter.
Let $\sigma$ be a $\zeta$-representation of $N$;
for $\gamma \in \Gamma$, define $\sigma^\gamma$ as in Lemma \ref{lemma6}.

(i) 
The mapping
\begin{equation*}
\sigma^\gamma : N \longrightarrow \mathcal U (\mathcal H_{\sigma})
\end{equation*}
is a  $\zeta$-representation of $N$.

(ii) We have
\begin{equation}
\label{eq11}
% eq11
\sigma^{\gamma_1\gamma_2} \, = \, ( \sigma^{\gamma_1} )^{\gamma_2} 
\end{equation}
for all $\gamma_1, \gamma_2 \in \Gamma$.
\end{Prop}

\emph{Proof.} 
Claim (i) follows from Lemma \ref{lemma6},
because $\zeta$ is normalised,
and checking Claim (ii) is straightforward.
\hfill $\square$

\medskip

\noindent \textbf{Second standing assumption.} 
\emph{All multipliers appearing from now on in this paper
are assumed to be normalised.}

\medskip

% Equation (\ref{eq8}) and Proposition \ref{prop5}
% motivate the following

It is convenient to define now projective analogues of
$\Gamma$-kernels, 
and $\Gamma$-faithfulness, a notion already used
in the formulation of Theorem~\ref{Thm4}.

\medskip

\noindent \textbf{Definitions.} 
\emph{Let $\Gamma$ be a group, $N$ a normal subgroup,
$\zeta \in Z^2(\Gamma, \mathbf T)$ a multiplier,
and $\sigma : N \longrightarrow \mathcal U (\mathcal H)$
a $\zeta$-representation of $N$.}

(i) \emph{The} projective $\Gamma$-kernel \emph{of $\sigma$ 
is the normal subgroup
\begin{equation}
\label{eq12}
% eq12
\aligned
\operatorname{P\ker_{\Gamma}}(\sigma) \, &= \,
\left\{
x \in \operatorname{Pker}(\sigma) 
\hskip.1cm \big\vert \hskip.1cm
\sigma^\gamma(x) = \sigma(x)
\hskip.1cm \text{for all} \hskip.1cm
\gamma \in \Gamma 
\right\} 
\\
&= \, \operatorname{Pker} \Big( \bigoplus_{\gamma \in \Gamma} \sigma^\gamma \Big) .
\endaligned
\end{equation}
of $\Gamma$.}

(ii) \emph{The projective representation $\sigma$ is}
$\Gamma$-P-faithful \emph{if $\operatorname{P\ker_{\Gamma}}(\sigma) \, = \, \{e\}$.}

\medskip

\noindent \textbf{Remarks.}
(A) 
In the particular case $\zeta = 1$, observe that
$\operatorname{P\ker_{\Gamma}}(\sigma)$
is a subgroup of $\operatorname{\ker_{\Gamma}}(\sigma)$
which can be a proper subgroup. 

(B)
Suppose that $N$ is a central subgroup in $\Gamma$.
For a $\zeta$-representation $\sigma$ of $N$, we have
\begin{equation*}
\operatorname{P\ker_{\Gamma}}(\sigma) \, = \,
\left\{
x \in \operatorname{Pker}(\sigma) 
\hskip.1cm \big\vert \hskip.1cm
\zeta(\gamma x,\gamma^{-1}) \zeta(\gamma,x) = 1
\hskip.1cm \text{for all} \hskip.1cm
\gamma \in \Gamma 
\right\} .
\end{equation*}
Since
\begin{equation*}
\zeta(\gamma x, \gamma^{-1})\zeta(x,\gamma) \, = \,
\zeta(x,\gamma) \zeta(x\gamma, \gamma^{-1})  \, = \,
\zeta(x, 1) \zeta(\gamma, \gamma^{-1}) \, = \, 1
\end{equation*} 
for every $x \in Z(\Gamma)$ and $\gamma \in \Gamma$
(recall that $\zeta$ is normalised), 
we have $\zeta(\gamma x,\gamma^{-1}) = \zeta(x,\gamma)^{-1}$
and therefore also
\begin{equation}
\label{eq13}
% eq13
\operatorname{P\ker_{\Gamma}}(\sigma) \, = \,
\left\{
x \in \operatorname{Pker}(\sigma) 
\hskip.1cm \big\vert \hskip.1cm
\zeta(x,\gamma) = \zeta(\gamma,x)
\hskip.1cm \text{for all} \hskip.1cm
\gamma \in \Gamma 
\right\} .
\end{equation}

\section{\textbf{ 
Extensions of  groups associated to multipliers
\\
Proof of Theorem \ref{Thm1}}}
\label{section4}
% s4

Consider a group $\Gamma$, a multiplier $\zeta \in Z^2(\Gamma, \mathbf T)$,
and a subgroup $A$ of $\mathbf T$ containing $\zeta(\Gamma \times \Gamma)$. 

We define a group $\Gamma(\zeta)$ with underlying set $A \times \Gamma$
and multiplication
\begin{equation}
\label{eq14}
% eq14
(s,x) (t,y) \, = \, (st\zeta(x,y), xy)
\end{equation}
for all $s,t \in A$ and $x,y \in \Gamma$;
observe that $(s,x)^{-1} = (s^{-1},x^{-1})$,
because $\zeta$ is normalised.
This fits naturally in a central extension
\begin{equation}
\label{eq15}
% eq15
\{e\} 
\hskip.2cm \longrightarrow \hskip.2cm
A \hskip.2cm \overset{s \mapsto (s,e)}{\longrightarrow} \hskip.2cm
\Gamma(\zeta) \hskip.2cm\overset{(s,x) \mapsto x}{\longrightarrow} \hskip.2cm
\Gamma \hskip.2cm \longrightarrow \hskip.2cm 
\{e\} .
\end{equation}
\emph{We insist on the fact 
that $\Gamma(\zeta)$ depends on the choice of $A$,
even if the notation does not show it.}
Whenever $H$ is a subgroup of $\Gamma$,
we identify $H(\zeta)$ with the appropriate subgroup of $\Gamma(\zeta)$.

\par

To any $\zeta$-representation $\pi$ of $\Gamma$ 
on some Hilbert space $\mathcal H$
corresponds a representation $\pi^0$ of $\Gamma(\zeta)$
on the same space
defined by
\begin{equation}
\label{eq16}
% eq16
\pi^0(s,x) \, = \, s\pi(x)
\end{equation}
for all $(s,x) \in \Gamma(\zeta)$.
Conversely, to any representation $\pi^0$ of $\Gamma(\zeta)$
on $\mathcal H$ which is the identity on $A$
(namely which is such that $\pi^0(a) = a \operatorname{id}_{\mathcal H}$
for all $a \in A$)
corresponds a $\zeta$-representation $\pi$ of $\Gamma$
defined by
\begin{equation}
\label{eq17}
% eq17
\pi(x) \, = \, \pi^0(1,x) .
\end{equation}

\begin{Lemma}
\label{lemma9}
% lem9
(i) 
The correspondance $\pi \leftrightsquigarrow \pi^0$
given by Equations (\ref{eq16}) and (\ref{eq17})
is a bijection between $\zeta$-representations of $\Gamma$ on $\mathcal H$
and representations of $\Gamma(\zeta)$ on $\mathcal H$ which are the identity on $A$.
\par

(ii) 
If $\pi$ and $\pi^0$ correspond to each other in this way, 
$\pi$ is irreducible if and only if $\pi^0$ is so.
\par

Assume moreover that the subgroup $A$ of $\mathbf T$ contains both
the image $\zeta(\Gamma \times \Gamma)$ of the multiplier and the subset
\begin{equation}
\label{eq18}
% eq18
T_{\pi} \, :=  \, 
\left\{ 
z \in \mathbf T 
\hskip.2cm \vert \hskip.2cm
z \operatorname{id}_{\mathcal H} = \pi(x)
\hskip.1cm \text{for some}  \hskip.1cm 
x \in \operatorname{Pker}(\pi)
\right\} 
\end{equation}
of $\mathbf T$.
\par

(iii) 
If $\pi$ and $\pi^0$ are as above, $\pi$ is  P-faithful if and only if $\pi^0$ is faithful.
\end{Lemma}

\emph{Observation.}
If $\Gamma$ is countable, 
$\zeta(\Gamma \times \Gamma)$ and $T_{\pi}$ 
are countable subsets of $\mathbf T$,
so that there exists a countable group $A$ as in (\ref{eq15})
which contains both $T_{\pi}$ and the image of  $\zeta$.

\medskip

\emph{Proof.}
Claims (i) and (ii) are obvious. 
The generalisation of Claim (i) for continuous representations
of locally compact groups appears as a corollary to Theorem~1 in \cite{Klep--74};
see also Theorem 2.1  in  \cite{Mack--58}.

For Claim (iii), suppose first that $\pi$ is P-faithful.
If $(s,x) \in \ker (\pi^0)$, namely if $s\pi(x) = 1$,
then $x \in \operatorname{Pker}(\pi)$, so that $x = e$;
it follows  that $s=1$, so that $(s,x) = (1,e)$.
Thus $\pi^0$ is faithful.
\par

Suppose now that $\pi^0$ is faithful.
If $x \in \operatorname{Pker}(\pi)$, namely if there exists $s \in \mathbf T$
such that $s\pi(x) = 1$, then $(s,x) \in \ker (\pi^0)$,
so that $s=1$ and $x=e$.
Thus $\pi$ is P-faithful.
\hfill $\square$

\medskip

\noindent
\textbf{Proof of Theorem \ref{Thm1}.}
Let $\pi$ be an irreducible P-faithful projective representation of $\Gamma$,
of multiplier $\zeta$.
Choose a subgroup $A$ of $\mathbf T$ 
containing $\zeta(\Gamma \times \Gamma)$ and $T_{\pi}$
(as defined in Lemma \ref{lemma9}).
Let $\Gamma(\zeta)$ be as in (\ref{eq14})
and $\pi^0$ be as in (\ref{eq16}).
Since $\pi^0$ is irreducible and faithful (Lemma \ref{lemma9}),
Schur's Lemma implies that 
$A$ is the centre  of $\Gamma(\zeta)$,
so that $\Gamma \approx  \Gamma(\zeta) / Z(\Gamma(\zeta))$.
\par

If $\Gamma$ is countable, 
$\Gamma(\zeta)$ can be chosen countable,
by the observation just after Lemma \ref{lemma9}.
\par

Conversely, let $\Delta$ be a group such that $\Gamma \approx \Delta / Z(\Delta)$
and let $\pi^0$ be a representation of $\Delta$ which is irreducible and faithful.
Again by Schur's Lemma, the subgroup 
$(\pi^0)^{-1}(\mathbf T )$ coincides with the centre of $\Delta$.
Let $\mu : \Gamma \longrightarrow \Delta$ be any set-theoretical section of the projection
$\Delta \longrightarrow \Delta / Z(\Delta) \approx \Gamma$, 
with $\mu(e_{\Gamma}) = e_{\Delta}$.
The assignment $\pi : \gamma \longmapsto \pi^0(\mu(\gamma))$
defines a projective representation of $\Gamma$ which is irreducible and P-faithful,
by Lemma \ref{lemma9}.
\hfill $\square$

\medskip

Our next lemma reduces essentially to Lemma \ref{lemma9}.iii
if $N = \Gamma$. It  will be used in 
the proof of Lemma \ref{lemma13}.

\begin{Lemma}
\label{lemma10}
% lemma 10
Consider a normal subgroup $N$ of $\Gamma$
and a $\zeta$-representation $\sigma$ of $N$
in some Hilbert space $\mathcal H$.
Let $A$ be a subgroup of $\mathbf T$ containing both
$\zeta(N \times N)$ and the subset
\begin{equation}
\label{eq19}
% eq19
T_{\sigma, \Gamma} \, :=  \, 
\left\{ 
z \in \mathbf T 
\hskip.2cm \vert \hskip.2cm
z \operatorname{id}_{\mathcal H} = \sigma(x)
\hskip.1cm \text{for some}  \hskip.1cm 
x \in \operatorname{Pker}_{\Gamma}(\sigma)
\right\} 
\end{equation}
of $\mathbf T$ (compare with Equation (\ref{eq18})).
Define $N(\zeta)$ and $\Gamma(\zeta)$ as in the beginning of the present section.
Then
\begin{itemize}
\item[(i)]
$(\sigma^{\gamma})^0 = (\sigma^0)^{\gamma}$ for all $\gamma \in \Gamma$ ;
\item[(ii)]
$\operatorname{P\ker_{\Gamma}}(\sigma) =
\left\{ x \in N \hskip.1cm \Big\vert \hskip.1cm
\aligned
&\text{there exists} \hskip.2cm s \in A 
\\
&\text{with} \hskip.2cm
(s,x) \in \ker_{\Gamma(\zeta)} (\sigma^0) 
\endaligned \right\}$,
so that, in particular, $\sigma$ is $\Gamma$-P-faithful if and only if
$\sigma^0$ is $\Gamma(\zeta)$-faithful.
\end{itemize}
\end{Lemma}

\emph{Proof.}
Checking (i) is straightforward.
\par

To show (ii), let  $x \in \operatorname{P\ker_{\Gamma}} (\sigma)$.
Thus there exists $s \in A$ such that
$\sigma^{\gamma} (x)= s^{-1} \operatorname{id}_{\mathcal H}$ 
for all $\gamma \in \Gamma$.
Then, for all $\gamma \in \Gamma$, we have
\begin{equation*}
(\sigma^0)^\gamma (s,x)
\, =\, (\sigma^\gamma)^0 (s,x)
\, = \, s\sigma^\gamma(x)
\, = \, s s^{-1}\operatorname{id}_{\mathcal H}
\,  = \, \operatorname{id}_{\mathcal H} ,
\end{equation*}
that is,  $(s,x)\in \ker_{\Gamma(\zeta)} (\sigma^0)$.
\par

Conversely, let $x \in N$  be such that there exists 
$s \in A$ with $(s,x)\in \ker_{\Gamma(\zeta)} \sigma^0$.
Then, for all $\gamma \in \Gamma$, we have
\begin{equation*}
\sigma^{\gamma}(x) 
\, = \,
(\sigma^{\gamma})^0 (1,x)
\, = \,
s^{-1} (\sigma^{\gamma})^0 (s,x)
= s^{-1} (\sigma^0)^{\gamma} (s,x)
\, = \, 
s^{-1} \operatorname{id}_{\mathcal H} ,
\end{equation*}
that is, $x\in \operatorname{P\ker_{\Gamma}} (\sigma)$.
\hfill $\square$

\medskip

Given a group $\Gamma$ and a multiplier $\zeta \in Z^2(\Gamma, \mathbf T)$,
a \emph{$\zeta$-character} of $\Gamma$ is a $\zeta$-representation 
$\chi : \Gamma \longrightarrow \mathbf T = \mathcal U (\mathbf C)$;
we denote by $X^{\zeta}(\Gamma)$ the set of all these.
Observe that, for $\chi_1, \chi_2 \in X^{\zeta}(\Gamma)$,
the product $\chi_1 \overline{\chi_2}$ is a character of $\Gamma$
in the usual sense,
namely a homomorphism from $\Gamma$ to $\mathbf T$.
Such a homomorphism factors via the \emph{abelianisation}
$\Gamma / [\Gamma,\Gamma]$, that we denote by
$\Gamma^{\operatorname{ab}}$.
We denote by $\widehat{\Gamma^{\operatorname{ab}}}$ the
\emph{character group} $\operatorname{Hom}(\Gamma, \mathbf T)$.
For further reference, 
we state here the following straightforward observations.

\begin{Lemma}
\label{lemma11}
% lemma 11
Let $\Gamma$ be a group and let $\zeta \in Z^2(\Gamma,\mathbf T)$.
\par
(i) If $\zeta \notin B^2(\Gamma, \mathbf T)$, then
$X^{\zeta}(\Gamma) = \emptyset$.
\par

(ii) If  $\zeta \in B^2(\Gamma, \mathbf T)$,
there exists a bijection between $X^{\zeta}(\Gamma)$ 
and $\widehat{\Gamma^{\operatorname{ab}}}$.

\iffalse
(iii) In particular, if $\zeta \in B^2(\Gamma, \mathbf T)$
and $\Gamma^{\operatorname{ab}}$ is finite,
$X^{\zeta}(\Gamma)$
is of the same order as $\Gamma^{\operatorname{ab}}$.
\fi
\end{Lemma}

\emph{Proof.}
(i) Suppose that there exists $\chi \in X^{\zeta}(\Gamma)$;
then $\zeta(x,y) = \frac{ \chi(x) \chi(y) }{ \chi(xy) }$ for all $x,y \in \Gamma$,
so that $\zeta$ is a coboundary.
\par

(ii) If $\zeta \in B^2(\Gamma, \mathbf T)$, there exists
a mapping $\nu : \Gamma \longrightarrow \mathbf T$
such that $\zeta$ is related to $\nu$ as in Fromula (\ref{eq5}),
so that $\nu \in X^{\zeta}(\Gamma)$.
For any $\chi \in X^{\zeta}(\Gamma)$, observe that
$\chi \overline{\nu}$ is an ordinary character of $\Gamma$,
so that $\chi \longmapsto \chi \overline{\nu}$
is a bijection 
$X^{\zeta}(\Gamma) \longrightarrow X^{1}(\Gamma) = 
\widehat{\Gamma^{\operatorname{ab}}}$.
\hfill $\square$

\section{\textbf{
Proof of  
$(i) \Longrightarrow (ii) \Longleftrightarrow (iii)$ in 
Theorem \ref{Thm4}}}
\label{section5}
% s5

The first proposition of this section is a reminder of
Section 3 of \cite{Mack--58}.

\begin{Prop}[Mackey]
\label{PropMackey}
A $\zeta$-representation $\pi$ of a countable group $\Gamma$ has
a direct integral decomposition in irreducible $\zeta$-representations, 
of the form
\begin{equation}
\label{eq20}
% eq20
\pi \, = \, \int_{\Omega}^{\oplus} \pi_{\omega} \, d\mu(\omega) .
\end{equation}
\end{Prop}

\emph{Proof.}
Consider  a subgroup $A$ of $\mathbf T$ 
containing $\zeta(\Gamma \times \Gamma)$
and
% \footnote{Mel Bachir 17 juin : il n'est pas n\'ecessaire que $A$
% contienne $T_\pi$.}
the subset $T_{\pi}$ defined in (\ref{eq18}),
the resulting extension $\Gamma(\zeta)$,
and the representation $\pi^0$ of $\Gamma(\zeta)$ defined in (\ref{eq16}).
There exists a direct integral decomposition in irreducible representations
\begin{equation*}
\pi^0 \, = \, \int_{\Omega}^{\oplus} (\pi^0)_{\omega} \, d\mu(\omega) 
\end{equation*}
with respect to a measurable field $\omega \longmapsto (\pi^0)_{\omega}$
of irreducible representations of $\Gamma(\zeta)$
on a measure space $(\Omega, \mu)$;
see \cite{Di--69C$^*$}, Sections 8.5 and 18.7.
\par

Since $\pi^0(s,x) = s\pi(x)$ for all $(s,x) \in \Gamma(\zeta)$, 
we have $(\pi^0)_{\omega}(s,x) = s \, (\pi^0)_{\omega}(1,x)$
for all $(s,x) \in \Gamma(\zeta)$ and for almost all $\omega \in \Omega$.
It follows that, for almost all $\omega \in \Omega$,  
the mapping
$\pi_{\omega} : \Gamma \longrightarrow \mathcal U (\mathcal H_{\omega})$
defined by $\pi_{\omega}(x) = (\pi^0)_{\omega}(1,x)$
is a $\zeta$-representation of $\Gamma$ which is irreducible, 
and $(\pi_{\omega})^0 = (\pi^0)_{\omega}$.
Hence we have a decomposition as in (\ref{eq20}).
\hfill $\square$

\medskip

We isolate in the next lemma an argument that we will use
in the proofs of Propositions \ref{Prop14}, \ref{Prop15}, and \ref{PropSection5}.

\medskip

\noindent
\textbf{Notation.} Let $\Gamma$ be a group and $N$ a normal subgroup.
We denote by $(C_j)_{j \in J}$ the $\Gamma$-conjugacy classes contained in $N$
and, for each $j \in J$, by $N_j$ the normal subgroup of $\Gamma$
generated by $C_j$.

\begin{Lemma}
\label{lemma13}
% lemma 13
 Let $\Gamma$ be a  group, $\zeta \in Z^2(\Gamma, \mathbf{T})$,
and $N$ a normal subgroup of $\Gamma$.
Let $(C_j)_{j\in J}$  and $(N_j)_{j \in J}$ be as above.
% Let $(C_j)_{j\in J}$  be an enumeration 
% of the $\Gamma$-conjucacy classes contained in $N$; 
% for each $ j\in J$ let $N_j$ denote 
% the normal subgroup of $\Gamma$ generated by $C_j$. 
Let $A$ be a subgroup of $\mathbf{T}$
containing $\zeta(\Gamma\times\Gamma)$, as well as $\chi(x)$
for every $\chi \in X^\zeta(N_j)$, $j \in J$, and $x \in N_j$.
Let $N(\zeta)$ be the  central extension of $N$
corresponding to $\zeta$ and $A$ as in (\ref{eq15}).
\par

Then, for every  $\zeta$-representation $\sigma$ of $N$, we have:
$\sigma$ is 
$\Gamma$-P-faithful if and only if the 
corresponding representation $\sigma^0$ of  $N(\zeta)$
is $\Gamma(\zeta)$-faithful.
\end{Lemma}

\emph{Remark.}
This lemma will be applied in situations where $\Gamma$ is a countable group.
Observe that, if $\Gamma$ is countable,
there exists a countable group $A$ as in the previous lemma
as soon as $\Gamma$ has Property (Fab),
or more generally as soon as $N_j^{\operatorname{ab}}$ is finite
for all $j \in J$ such that the restriction of $\zeta$ to $N_j$ 
is in $B^2(N_j, \mathbf T)$.

\medskip

\emph{Proof.} 
In view of Lemma \ref{lemma10}, it suffices to prove that
$A$ contains $T_{\sigma, \Gamma}$ 
for every  $\zeta$-representation $\sigma$ of $N$.
\par  

Let $z \in T_{\sigma, \Gamma}$; choose 
$x \in \operatorname{Pker}_{\Gamma} (\sigma)$ 
such that $\sigma(x) = \operatorname{id}_{\mathcal H_\sigma}$.
Let $j \in J$ be such that $x \in C_j$;
we have  $N_j \subset  \operatorname{Pker}_{\Gamma} (\sigma)$,
because the latter group  is  normal in $\Gamma$. 
The restriction of $\sigma$ to $N_j$ defines a $\zeta$-character
$\chi \in X^{\zeta}(N_j)$ such that $z = \chi(x)$. 
This shows that  $z \in A$, by the choice of $A$.
\hfill $\square$

\medskip

Implications  $(i) \Longrightarrow (ii)$ and $(i) \Longrightarrow (iii)$
of Theorem \ref{Thm4}
are particular cases of the following proposition, because
the minisocle $MS(\Gamma)$ and the subgroup $MA(\Gamma)$
of a countable group $\Gamma$
have the properties assumed for the group $N$ below
(Proposition 1 in  \cite{BeHa--08}).

\begin{Prop}
\label{Prop14}
% Proposition 14
Let $\Gamma$ be a countable group, $N$ a normal subgroup,
and $\zeta \in Z^2(\Gamma, \mathbf T)$.
Let $(C_j)_{j\in J}$  and $(N_j)_{j \in J}$ be as just before  Lemma \ref{lemma13}.
% Let $\big( C_j \big)_{j \in J}$ be an enumeration
% of the $\Gamma$-conjugacy classes contained in~$N$;
% for each $j \in J$ let $N_j$ denote the normal subgroup of $\Gamma$
% generated by $C_j$.
Assume that the abelianised group $N_j^{\operatorname{ab}}$ is finite for all $j \in J$
such that the restriction to $N_j$ of $\zeta$ is in $B^2(N_j,\mathbf T)$.
\par

Let  $\pi$ be a $\zeta$-representation of $\Gamma$
and let
\begin{equation}
\label{eq21}
% eq21
\sigma \, :=  \, \pi \vert_N \, = \, 
\int_{\Omega}^{\oplus} \sigma_{\omega} \, d\mu(\omega) 
\end{equation}
be a direct integral decomposition of the restriction of $\pi$ to $N$
in irreducible $\zeta$-representations $\sigma_{\omega}$ of $N$.
\par

If $\pi$ is irreducible and P-faithful,
then $\sigma_{\omega}$ is $\Gamma$-P-faithful
for almost all $\omega \in \Omega$.
\end{Prop}

\emph{Proof.}
The strategy is to reduce the proof to the case of
ordinary representations and to use  Lemma 9 of \cite{BeHa--08}.
\par

By hypothesis and by Lemma \ref{lemma11},  $X^{\zeta}(N_j)$  
is finite (possibly empty)  for all $j\in J$. 
Since $J$ is countable, 
we can choose a  \emph{countable}  subgroup $A$ of $\mathbf T$
containing  $\zeta(\Gamma \times \Gamma)$,
as well as $\chi(x)$
for every $\chi \in X^\zeta(N_j)$, $j \in J$, and $x \in N_j$.
\par

Let $\Gamma(\zeta)$ and $N(\zeta)$  be as in  (\ref{eq14});   
let $\pi^0$ and $\sigma_{\omega}^0$ be the representations
of $\Gamma(\zeta)$ and $N(\zeta)$ corresponding to 
the $\zeta$-representations $\pi$
and $\sigma_\omega$, respectively.
Because $\pi$ is P-faithful, 
the subset $T_{\pi}$ defined in (\ref{eq18}) is reduced to $\{e\}$
and therefore $\pi^0$ is faithful (Lemma \ref{lemma9}).

Since (see the proof of Proposition  \ref{PropMackey})
\begin{equation*}
\sigma^0
\, = \, \pi^0 \vert_N
\, = \, \int_{\Omega}^{\oplus} \sigma_{\omega}^0 d\mu(\omega) ,
\end{equation*}
the representation $\sigma_\omega^0$ of $N(\zeta)$ is 
$\Gamma(\zeta)$-faithful for almost all $\omega$ (Lemma 9
of  [BeHa--08]).  Therefore, by Lemma \ref{lemma13}, 
$\sigma_\omega$ is $\Gamma$-P-faithful for almost all $\omega$.
\hfill $\square$

\medskip

The equivalence  $(ii) \Longleftrightarrow (iii)$ 
of Theorem \ref{Thm4}
is a particular case 
of the following Proposition.

\begin{Prop}
\label{Prop15}
% Proposition 15
Assume that the normal subgroup $N$ of $\Gamma$
is a direct product $B \times S$ of normal subgroups of $\Gamma$,
and that $S = \prod_{i \in I} S_i$ is a restricted direct product
of finite simple nonabelian  subgroups $S_i$.
%  which are normal in $S$.
Assume moreover that any $\Gamma$-invariant subgroup of $B$
generated by one $\Gamma$-conjugacy class has finite abelianisation.
\par
The following conditions are equivalent:
\begin{itemize}
\item[($\alpha$)]
$N$ has a $\Gamma$-P-faithful irreducible $\zeta$-representation;
\item[($\beta$)]
$B$ has a $\Gamma$-P-faithful irreducible $\zeta$-representation.
\end{itemize}
\end{Prop}

\emph{Proof}
The proof of the implication $(\alpha) \Rightarrow (\beta)$
follows closely the  proof of  Proposition \ref{Prop14}, 
with one difference: one has to use
the more general version of Lemma 9 in \cite{BeHa--08} 
which is mentioned at the bottom of page 866 of this article.
\par

For the converse implication, we assume now that
$B$ has a $\Gamma$-P-faithful irreducible $\zeta$-representation $\sigma$.
The group $S$ has a faithful irreducible (unitary) representation, 
say $\rho$, such  that $\rho(x) \notin \mathbf T$ for all $x \in S, x \ne e$,
namely $\rho$ is P-faithful;
see the proof  of Lemma 13 in \cite{BeHa--08}
(this Lemma 13 contains a hypothesis "$A$ abelian",
but it is redundant for the part of the proof we need here).
The tensor product $\sigma \otimes \rho$
is an irreducible $\zeta$-representation of $N$.
Since $\sigma$ is $\Gamma$-P-faithful, it follows
from Lemma 12 of \cite{BeHa--08}
that $\sigma \otimes \rho$ is $\Gamma$-P-faithful.
\hfill $\square$

\section{\textbf{
End of proof of Theorem \ref{Thm4}}}
\label{section6}
% s6

Let us first recall the definition of induction for projective representations,
from Section 4 in \cite{Mack--58}.

Let $\Gamma$ be a group, 
$\zeta \in Z^2(\Gamma, \mathbf T)$ a multiplier,
$H$ a subgroup of $\Gamma$,
and $\sigma : H \longrightarrow \mathcal U (\mathcal K)$ a $\zeta$-representation.
Let $\mathcal H$ be the Hilbert space of mappings
$f : \Gamma \longrightarrow \mathcal K$ with\footnote{We use
$H \backslash \Gamma$, rather than $\Gamma / H$ as in \cite{BeHa--08},
which provides easier formulas.}
% gauche$\leftrightsquigarrow$droite.}
the two following properties:
\begin{itemize}
\item[$\circ$] 
$f(hx) \, = \, \zeta(h,x) \, \sigma(h)(f(x))$ for all
$x \in \Gamma$ and $h \in H$,
\item[$\circ$] 
$\sum_{x \in \Gamma \backslash H} \Vert f(x) \Vert^2 \, < \, \infty$.
\end{itemize}
The $\zeta$-representation $\operatorname{Ind}_H^{\Gamma}(\sigma)$
of $\Gamma$ is the multiplier representation of $\Gamma$ in $\mathcal H$
defined by
\begin{equation}
\label{eq22}
% eq22
\left( \operatorname{Ind}_H^{\Gamma}(\sigma) \, (x) \, f \right) )(y) \, = \,
f(yx)
\end{equation}
for all $x,y \in \Gamma$.

It can be checked (see \cite{Mack--58}, Pages 273-4)
that the representation $\left( \operatorname{Ind}_H^{\Gamma}(\sigma) \right)^0$
of $\Gamma(\zeta)$ associated to $\operatorname{Ind}_H^{\Gamma}(\sigma)$
is the representation $\operatorname{Ind}_{H(\zeta)}^{\Gamma(\zeta)}(\sigma^0)$
induced by the representation $\sigma^0$ 
from  $H(\zeta)$ to $\Gamma(\zeta)$.

\medskip

The last claim of Theorem \ref{Thm4} follows from the next proposition.

\begin{Prop}
\label{PropSection5}
% Proposition 16
Let $\Gamma$ be a countable group and let $\zeta \in Z^2(\Gamma,\mathbf T)$.
Let $(C_j)_{j\in J}$  and $(N_j)_{j \in J}$ be as just before  Lemma \ref{lemma13},
with $N = \Gamma$.
% Let $\big( C_j \big)_{j \in J}$ be an enumeration
% of the conjugacy classes  in $\Gamma$;
% for each $j \in J$ let $N_j$ denote the normal subgroup of $\Gamma$
% generated by $C_j$.
Assume that the abelianised group $N_j^{\operatorname{ab}}$ is finite for all $j \in J$
such that the restriction to $N_j$ of $\zeta$ is in $B^2(N_j,\mathbf T)$.\footnote{This
assumption holds whenever $\Gamma$ has Property (Fab).}
\par

Let $\sigma$ be a $\zeta$-representation of the minisocle $MS(\Gamma)$.
Set $\pi := \operatorname{ind}^{\Gamma}_{MS(\Gamma)} (\sigma)$
and let
\begin{equation}
\pi \, = \, \int_{\Omega}^{\oplus} \pi_{\omega} d\mu (\omega)
\end{equation}
be a direct integral decomposition of $\pi$ in irreducible $\zeta$-representations
of $\Gamma$. 
\par

If $\sigma$ is  irreducible and $\Gamma$-P-faithful,
then $\pi_{\omega}$ is P-faithful for almost all $\omega \in \Omega$.
\end{Prop}

\emph{Proof.}
As for Proposition \ref{Prop14},
the strategy is to reduce the proof 
to the case of ordinary representations, and to use this time
Lemma 10 of \cite{BeHa--08}.
We write $M$ for $MS(\Gamma)$.

\par

By hypothesis and by Lemma \ref{lemma11}, 
we can choose a countable subgroup $A$ be a $\mathbf T$
containing the sets $\zeta(\Gamma \times \Gamma)$  and
$X^{\zeta}(N_j)(x)$ for every $j \in J$ and every $x \in N_j$.
We consider the corresponding  extension $\Gamma(\zeta)$ of $\Gamma$. 
Denote by $\pi^0$ and $\pi^0_{\omega}$ the
representations of $\Gamma(\zeta)$ corresponding to the
$\zeta$-representations $\pi$ and $\pi_{\omega}$
of $\Gamma$, and similarly $\sigma^0$
for the representation  of $M(\zeta)$ corresponding to the
$\zeta$-representation $\sigma$ of $M$.

We have
\begin{equation*}
\pi^0
\, = \, {\rm ind}_{M(\zeta)}^{\Gamma (\zeta)} \sigma^0 
\, = \,  \int_{\Omega}^{\oplus} \pi_\omega^0 d\mu(\omega).
\end{equation*}
In view of  Lemma \ref{lemma13} applied to $N=\Gamma$,  
it suffices to show that $\pi_{\omega}^0$ is P-faithful for almost all $\omega$.

Since $\sigma$ is $\Gamma$-P-faithful, 
we have that $\sigma^0$ is $\Gamma$-faithful
(again by Lemma \ref{lemma13}).  
It will follow from Lemma 10 in \cite{BeHa--08}
that $\pi_{\omega}^0$ is faithful for almost all $\omega$
provided we show that  $M(\zeta) \cap L \neq  \{ e\}$, 
for every finite foot $L$ in $\Gamma(\zeta)$. 

In order to check this condition, 
let $L$ be  finite foot  in $\Gamma(\zeta)$. 
We claim that $L\subset M(\zeta)$. 
Indeed, recall that, set-theoretically, 
we have $\Gamma(\zeta)= A \times \Gamma$ and  $M(\zeta)= A \times M$; 
thus, for any $(t,y) \in L$
with $(t,y) \neq e$, we have $y \in M (= MS(\Gamma))$, 
and therefore $(t,y )\in M(\zeta)$.
\hfill $\square$

\section{\textbf{
Capable and incapable groups
\\
Proof of Proposition \ref{PropBFSE}}}
\label{section7}

In  Proposition \ref{PropBFSE}, 
Claims~(i) and (ii) are respectively Corollary 2.3 
and part of Corollary 2.2 of \cite{BeFS--79}.
% is Proposition 3.9, Page 209 of \cite{BeTa--82}
Claim (iii) is a consequence of Claim (i), 
in a formulation and with a proof 
shown to us by Graham Ellis \cite{Ellis}, see below.
Corollary \ref{Cor3} is an immediate consequence of
Theorem \ref{Thm1} and Proposition \ref{PropBFSE}.

\medskip

\noindent
\textbf{Proof of Claim (iii) in Proposition \ref{PropBFSE}.}
For a central extension
\begin{equation*}
\{e\} \hskip.1cm \to \hskip.1cm
A \hskip.1cm \to \hskip.1cm
\Gamma \hskip.1cm \to \hskip.1cm
\Gamma/A \hskip.1cm \to \hskip.1cm
\{e\} ,
\end{equation*}
the Ganea extension of the Hochschild-Serre exact sequence 
in homology with trivial coefficients $\mathbf Z$ is
\begin{equation*}
\setcounter{MaxMatrixCols}{20}
\begin{matrix}
&& A \otimes_{\mathbf Z} \Gamma^{\operatorname{ab}}
&& \longrightarrow && H_2(\Gamma,\mathbf Z)
&& \longrightarrow &&  H_2(\Gamma/A,\mathbf Z)
&& \longrightarrow  &
\\
&& A
&& \longrightarrow  && H_1(\Gamma,\mathbf Z)
&&  \longrightarrow  && H_1(\Gamma/A,\mathbf Z)
&& \longrightarrow &\{0\} 
\end{matrix}
\end{equation*}
(see for example \cite{EcHS--72}).
If $\Gamma$ is perfect (so that $\Gamma^{\operatorname{ab}} = \{0\}$), 
this reduces to
\begin{equation*}
\{0\} \hskip.2cm \longrightarrow \hskip.2cm
H_2(\Gamma, \mathbf Z) \hskip.2cm \longrightarrow \hskip.2cm
H_2(\Gamma/A, \mathbf Z) \hskip.2cm \longrightarrow \hskip.2cm
A \hskip.2cm \longrightarrow \hskip.2cm \{0\}
\end{equation*}
and it follows from the definition of the epicentre of $\Gamma$
that $Z^*(\Gamma) = Z(\Gamma)$.
Thus Claim (iii) is a straightforward consequence of Claim (i)
of Proposition \ref{PropBFSE}.
\hfill $\square$

\medskip

It is well-known that any cyclic group $C \ne \{e\}$ is incapable.
Indeed, suppose \emph{ab absurdo} that $C = \Delta / Z(\Delta)$.
Choose a generator $s$ of $C$ and a preimage $t$ of $s$ in $\Delta$;
any $\delta \in \Delta$ can be written as $\delta = z t^j$
for some $z \in Z(\Delta)$ and $j \in \mathbf Z$,
and two elements of this kind commute with each other,
so that $\Delta$ is abelian, hence $Z(\Delta) = \Delta$,
incompatible with $C \ne \{e\}$. 
The next lemma, which appears on Page 137 of \cite{Hall--40},
rests on an elaboration of the same argument.

\begin{Lemma}
\label{lemmaHall}
Let $\Gamma$ be a group containing
an element $s_0 \ne e$ such that the set
\begin{equation*}
\{ s \in \Gamma 
\hskip.2cm \vert  \hskip.2cm 
\text{there exists $n \in \mathbf Z$ with $s^n = s_0$} \}
\end{equation*}
(where $n$ can depend on $s$)
generates $\Gamma$.
Then $\Gamma$ is incapable.
\end{Lemma}

\emph{Proof.}
It suffices to show that, given any central extension
\begin{equation*}
\{e\} \hskip.2cm \longrightarrow \hskip.2cm
A \hskip.2cm \longrightarrow \hskip.2cm
\Delta \hskip.2cm \overset{\pi}{\longrightarrow} \hskip.2cm
\Gamma \hskip.2cm \longrightarrow \hskip.2cm 
\{e\} ,
\end{equation*}
$s_0$ has a preimage $t_0$ in $\Delta$ which  is central.
\par

Let $\delta \in \Delta$.
There exists 
\begin{equation*}
s_1,\hdots,s_k \in \Gamma
\hskip.5cm \text{and} \hskip.5cm
j_1,\hdots,j_k,n_1,\hdots,n_k \in \mathbf  Z
\end{equation*}
such that 
\begin{equation*}
\pi(\delta) = s_1^{j_1} \cdots s_k^{j_k}
\hskip.5cm \text{and} \hskip.5cm
s_1^{n_1} = \cdots = s_k^{n_k} = s_0 .
\end{equation*}
For $i = 0, \hdots, k$, choose a preimage $t_i$ of $s_i$ in $\Delta$.
There exist $a,a_1,\hdots, a_k \in A$ with
\begin{equation*}
\delta = a  t_1^{j_1} \cdots t_k^{j_k}
\hskip.5cm \text{and} \hskip.5cm
a_i t_i^{n_i} = t_0
\hskip.2cm \text{for} \hskip.2cm
i = 1,\hdots, k .
\end{equation*}
It follows that $t_0$ commutes with $t_i$ for $i = 1, \hdots, k$,
and thus that $t_0$ commutes with $\delta$, as was to be shown.
\hfill $\square$

\medskip

Claims (i) to (v) of the following proposition
are straightforward consequences of Lemma \ref{lemmaHall},
and the last claim follows from Theorem \ref{Thm1}.

\begin{Prop}
\label{PropIncapables}
% Proposition 18
The following groups are incapable:
\begin{itemize}
\item[(i)]
cyclic groups,  quasicyclic groups $\mathbf  Z (p^{\infty})$,
and the groups $\mathbf Z \left[ 1/m \right]$ for all integers $m \ge 2$;
% [more generally, $\mathbf  Z (p^{\infty}) \times H$ is incapable
% for any $p$--group $H$ \cite{BeFeS--79, Page 163}];
% DEJA DIT !!!
\item[(ii)]
finite abelian groups 
$\mathbf  Z / d_1\mathbf  Z \times \cdots \times \mathbf  Z / d_m\mathbf  Z$,
\item[]
(where $n \ge 2$, $d_1,\hdots,d_m \ge 2$, $d_1 \vert d_2 \vert \cdots \vert d_m$)
\item[]
with $d_{m-1} < d_m$;
\item[(iii)]
subgroups of $\mathbf  Q$;
\item[(iv)]
$\operatorname{SL}_2(\mathbf  Z) = \langle s,t 
\hskip.2cm \vert \hskip.2cm
s^2 = t^3 
\hskip.2cm \text{is central of order} \hskip.2cm 2 \rangle$;
\item[(v)]
$\langle s,t 
\hskip.2cm \vert \hskip.2cm
s^m = t^n 
\hskip.2cm \text{and} \hskip.2cm (s^m)^k = 1 \rangle$ \hskip.1cm
for $m,n \ge 1$, $k \ge 2$, 
\item[]
as well as
$\langle s,t 
\hskip.2cm \vert \hskip.2cm
s^m = t^n  \rangle$.
% \item[(vi)]
% ................................................................... (???)
\end{itemize}

\par
In particular, these groups do not afford any
irreducible P-faithful projective representation.
\end{Prop}

\emph{Remarks.}
About (i):
for any prime $p$, the quasicyclic group  $\mathbf  Z (p^{\infty})$
is the subgroup of $\mathbf T$
of roots of $1$ of order some power of $p$;
equivalently, $\mathbf  Z (p^{\infty})$ is the quotient
$\mathbf Q_p / \mathbf Z_p$
of the $p$-adic numbers 
by the $p$-adic integers.
\par

About (iii): there is a classification 
of the subgroups of $\mathbf  Q$, which is is standard;
see for example Chapter 10 of \cite{Rotm--95}.
\par

About (iv), let us recall that $\operatorname{SL}_2(\mathbf Z)$ 
is generated by a square root
$
\left( \begin{matrix}
\phantom{-}0 & 1 \\ -1 & 0 
\end{matrix} \right)
$
and a cubic root
$
\left( \begin{matrix}
0 & -1 \\ 1 &  \phantom{-}1
\end{matrix} \right)
$
of the central matrix
$
\left( \begin{matrix}
-1 & \phantom{-}0 \\ \phantom{-}0 & -1 
\end{matrix} \right) .
$
Next, since 
$\operatorname{SL}_2(\mathbf  Z) = 
\langle s,t \hskip.2cm \vert \hskip.2cm
s^2 = t^3 ,  \hskip.1cm s^4 = 1 \rangle$
has deficiency $\ge 0$
and finite abelianisation (indeed 
$\operatorname{SL}_2(\mathbf  Z)^{\operatorname{ab}}
\approx \mathbf Z / 12\mathbf Z$), 
it follows from Philip Hall's Inequality\footnote{Namely:
for a finitely presented group $\Gamma$,
the deficiency of $\Gamma$ is bounded by the difference
$\dim_{\mathbf Q}((\Gamma)^{\operatorname{ab}} \otimes_{\mathbf Z} \mathbf Q)
- s(H_2(\Gamma, \mathbf Z))$,
where $s(H)$ stands for the minimum number of generators of the group $H$;
see e.g. Lemma 1.2 in \cite{Epst--61}.
Recall also that the \emph{deficiency} of a finite presentation of a group
is the number of its generators minus the number of its relations,
and the deficiency of a finitely presented group the maximum 
of the deficiencies of its finite presentations.}
% (Lemma 1.2 in \cite{Epst--61}) 
that
$H_2(\operatorname{SL}_2(\mathbf Z), \mathbf Z) = \{0\}$.
\par

Similarly, $H_2(\operatorname{PSL}_2(\mathbf Z), \mathbf Z) = \{0\}$.
This follows alternatively from the formula
$H_n(\Gamma_1 \ast \Gamma_2, \mathbf Z) \approx
H_n(\Gamma_1, \mathbf Z) \oplus H_n(\Gamma_2, \mathbf Z)$
for $n \ge 1$, see Corollary 6.2.10
% page 170
in \cite{Weib--94}.
But $\operatorname{PSL}_2(\mathbf Z)$ is capable, since
its centre is trivial.

About  (v): the group $\operatorname{SL}_2(\mathbf Z)$ is of course
a particular case of groups in (v);
if $m$ and $n$ are coprime and at least $2$, the group
$\langle s,t  \hskip.2cm \vert \hskip.2cm s^m = t^n  \rangle$
is a \emph{torus knot group}.
\par

\section{\textbf{
On abelian groups}}
\label{section8}

The next proposition rests on a construction which appears in many places,
including \cite{Mack--49} and \cite{Weil--64}.
% [Loomis], etc.
It is part of the \emph{Stone--von Neumann--Mackey Theorem},
see the beginning of  \cite{MuNN--91}.
\par

Let $L$ be an abelian group, written multiplicatively.
Consider the group $X(L) = \operatorname{Hom}(L,\mathbf T)$ of characters of $L$,
with the topology of the simple convergence, 
which makes it a locally compact abelian group.
By Pontryagin duality, we can (and do) identify $L$
to the group of continuous characters on $X(L)$.
Consider also a dense subgroup $M$ of $X(L)$
and the direct product group $L \times M$. 
The mapping
\begin{equation*}
\zeta \, : \, 
(L \times M) \times (L \times M) \longrightarrow \mathbf T ,
\hskip.2cm
((\ell,m),(\ell',m')) \longmapsto m'(\ell)
\end{equation*}
is a multiplier on $L \times M$. 
Let $A$ be a subgroup of $\mathbf T$ containing the image of $\zeta$.
By definition, the corresponding
\emph{generalised Heisenberg group} is
\begin{equation*}
H^A_{L,M} \, = \, A \times L \times M
\end{equation*} 
with product defined by
\begin{equation*}
(z,\ell,m) \, (z',\ell',m') \, = \,
(zz' \ell (m') \, , \, \ell \ell' \, , \, mm') .
\end{equation*}
It is routine to check that the centre of $H^A_{L,M}$
is $A$.

\begin{Prop}
\label{PropSvN}
% Proposition 19
Any abelian group of the form $L \times M$,
with $M$ dense in  $X(L)$ as above,
affords a projective representation with is irreducible and P-faithful.
\end{Prop}

\emph{Proof.}
Let us sketch the definition and some properties of the
``Stone--von Neumann--Mackey representation'' of $H^A_{L,M}$
on  $\ell^2(L)$;
the latter is a Hilbert space,
with scalar product defined by
$\langle \xi \mid \eta \rangle = \sum_{\ell \in L} \overline{\xi(\ell)} \eta(\ell)$.
\par

For $(z,\ell,m)  \in H^A_{L,M}$ and $\xi \in \ell^2(L)$, set
\begin{equation*}
(R(z,\ell,m)\xi)(x) \, = \, z  m(x) \xi(x\ell) 
\hskip.5cm \text{for all} \hskip.2cm  x \in L .
\end{equation*}
It can be checked that $R(z,\ell,m)$ is a unitary operator on ${\mathcal H}=\ell^2(L)$ 
and that
\begin{equation*}
R \, : \, H^A_{L,M} \longrightarrow \mathcal U (\mathcal H)
\end{equation*}
is a representation of $H^A_{L,M}$ on $\mathcal H$.
\par

The space ${\mathcal H}$ has a natural orthonormal basis
$\left( \delta_u \right)_{u \in L}$.
It is easy to check that
\begin{equation*}
R(z,\ell,m) \delta_u \, = \, z m(u\ell^{-1})\delta_{u\ell^{-1}} ,
\end{equation*}
so that the representation $R$ is faithful.
Observe that, for all $u \in L$ and $m \in M$, the vector $\delta_u$
is an eigenvector of $R(1,1,m)$ with eigenvalue $m(u)$. 
If
\begin{equation*}
V_u \, = \, \{ \xi \in {\mathcal H} \mid R(1,1,m)\xi = m (u)\xi
\hskip.2cm \text{for all} \hskip.2cm m \in M \} ,
\end{equation*}
then $V_u = \mathbf C \delta_u$ is an eigenspace of dimension $1$
and ${\mathcal H} = \bigoplus_{u \in L} V_u$ (Hilbert sum).
\par

Let now $S \in \mathcal L (\mathcal H)$ 
be an operator commuting with  $R(1,1,m)$ for all $m \in M$.
Since $M$ is dense in $X(L)$, for every $u,v\in L$ with $u\neq v$, there exists
$m\in M$  such that $m(u)\neq m(v)$. 
As is easily checked,  this implies  that $S$ is diagonal
with respect to the basis $\left( \delta_u \right)_{u \in L}$, 
namely that there exist complex numbers $s_u$
such that $S(\delta_u) = s_u \delta_u$ for all $u \in M$.
Suppose moreover that $S$ commutes with  $R(1,\ell,1)$ for all $\ell \in L$;
since $R(1,\ell,1) \delta_u = \delta_{u\ell^{-1}}$, we have
$s_u = s_{u\ell^{-1}}$ for all $u,\ell \in L$.
Thus $S$ is a scalar multiple of the identity operator.
It follows from Schur's lemma that the representation $R$ is irreducible.
\par 

The representation $R$ of $H^A_{L,M}$ provides
a projective representation  of $L \times M$
which is irreducible and P-faithful. 
\hfill $\square$

\medskip

In particular, the following groups afford
projective representations which are irreducible and P-faithful:
\begin{itemize}
\item[$\circ$]
$\mathbf Z^n$ for any $n \ge 2$, 
as $\mathbf Z^{n-1}$ is a dense subgroup of 
$X(\mathbf Z) \approx \mathbf T$.
\item[$\circ$]
$\mathbf Z (p^{\infty}) \times \mathbf Z$, 
as $\mathbf{Z}$ is a dense subgroup of
$X(\mathbf Z(p^{\infty})) \approx \mathbf Z_p$.
For the latter isomorphism, see e.g. \cite{Bour--67}, 
chap. 2, \S~1, no. 9, cor. 4 of prop. 12;
$\mathbf Z (p^{\infty})$ is as just after Proposition \ref{PropIncapables}.
\item[$\circ$] 
$\mathbf{Q}^n$ for any $n\geq 2$. 
Indeed, let us check this for $n=2$, 
the general case being entirely similar.
The group 
$X(\mathbf{Q})$ can be identified with $\mathbf{A}/ \varphi(\mathbf{Q})$;
here $\mathbf{A}$ is the group of adeles  of $\mathbf{Q}$
and $\varphi : \mathbf{Q} \longrightarrow \mathbf{A}$ 
is the diagonal  embedding of  $\mathbf{Q}$ in $\mathbf A$
(recall that $\varphi (\mathbf{Q})$ is discrete and cocompact in $\mathbf A$).
More precisely, let $\chi_0$ be a non-trivial character of  
$\mathbf{A}$ with $\chi_0\vert_{\varphi(\mathbf Q)}=1$. 
Then  the mapping  
\begin{equation*}
\Phi : \mathbf A \longrightarrow X(\mathbf{Q}), 
\qquad  
a \longmapsto \left(q \mapsto \chi_0(a\varphi(q)) \right)
\end{equation*} 
factorizes to an isomorphism
$\mathbf{A}/ \varphi(\mathbf{Q}) \to X(\mathbf{Q})$ 
(see Chapter 3 in  \cite{GelfandEtCo}).
Fix  $a_0 \in \mathbf{A}$ with $a_0 \notin  \varphi(\mathbf{Q})$ 
and define a group homomorphism 
\begin{equation*}
f: \mathbf{Q} \longrightarrow  X(\mathbf Q), 
\qquad 
f(q) \, = \,  \Phi (a_0\varphi(q)).
\end{equation*} 
Then $f$ is injective since $a_0\varphi(q)\notin  \varphi(\mathbf{Q})$
for all $q\in \mathbf Q^*$. 
We claim that the range of $f$ is dense. Indeed,  assume that this is not the case.
By Pontryagin duality,  there  exists
$ q_0\in \mathbf{Q}^*$ such that  $f(q)(\varphi(q_0))=1$ for all $q\in \mathbf Q$.
This means that $\chi_0(a_0 \varphi(q_0q))=1$ for all  $q\in \mathbf Q$,
that is,  $\Phi (a_0\varphi(q_0))$ is  the trivial character of $\mathbf Q$.
This is a contradiction, since $a_0 \varphi(q_0) \notin   \varphi(\mathbf{Q})$.
\end{itemize}
Note that Proposition \ref{PropSvN} carries over to dense subgroups
of groups of the form $B \times X(B)$, 
with $B$ a locally compact abelian group.

\medskip

The case of finite groups is covered by a result of Frucht \cite{Fruc--31}.
For a modern exposition (and improvements\footnote{Any
finite abelian group affords two irreducible projective representations
of which the direct sum is P-faithful.
For a characterisation of those finite groups which have
a faithful linear representation which is a direct sum 
of $k$ irreducible representations, see Page 245,
and indeed all of Chapter 9, in \cite{BeZh--98}.}), 
see Page 166 of \cite{BeZh--98}.

\begin{Prop}[Frucht]
\label{PropFrucht}
% Proposition 20
For a finite abelian group $\Gamma$, the two following properties are equivalent:
\begin{itemize}
\item[(i)]
$\Gamma$ affords a projective representation which is irreducible and P-faithful;
\item[(ii)]
there exists a (finite abelian) group $L$ such that $\Gamma$ is isomorphic to the
direct sum $L \times L$.
\end{itemize}
\end{Prop}

\emph{Observation}, from \cite{Sury--08}.
Consider a prime $p$, the ``Heisenberg group'' $H$ below,
and its noncyclic centre $Z(H)$:
\begin{equation*}
H \, = \, 
\left( \begin{matrix}
1 & \mathbf F_p & \mathbf F_{p^2} \\
0 & 1 & \mathbf F_{p^2} \\
0 & 0 & 1
\end{matrix} \right)
\, \supset \,
Z(H) \, = \, 
\left( \begin{matrix}
1 & 0 & \mathbf F_{p^2} \\
0 & 1 & 0 \\
0 & 0 & 1
\end{matrix} \right)
\, \approx \, 
\mathbf F_{p^2} 
\, \approx \, 
\mathbf F_p \oplus \mathbf F_p 
\end{equation*}
(where the last $\approx$ indicates of course an isomorphism of additive groups,
not of rings!).
The quotient  $H/Z(H) \approx \mathbf F_{p}  \oplus \mathbf F_{p^2}$
is an abelian group of which the order $p$ is not a square,
and therefore which does not have the properties of Proposition \ref{PropFrucht},
but which is however a capable group.

\medskip

Recall (from just after Proposition \ref{PropSvN}) that 
$\mathbf Z^3$ affords a projective representation which is irreducible and P-faithful,
and compare with Claim (ii) of Proposition \ref{PropFrucht}.

\section{\textbf{ 
Final remarks}}
\label{section9}

In some sense,
what follows goes back for finite groups  to papers by Schur, from 1904 and 1907.
For the general case, see Sections V.5 and V.6 in \cite{Stam--73}, 
\cite{EcHS--72}, and \cite{Kerv--70}.

A \emph{stem cover} of a group $\Gamma$ is a group $\widetilde \Gamma$
given with a surjection $p$ onto $\Gamma$ such that
$\ker (p)$ is central in $\widetilde \Gamma$, 
contained in $[\widetilde \Gamma, \widetilde \Gamma]$, 
and isomorphic to $H_2(\Gamma, \mathbf Z)$.
Any group has a stem cover. 
The isomorphism type of $\widetilde \Gamma$ is
uniquely determined  in case $\Gamma$ is perfect, but not in general.
For example, the dihedral group of order $8$ and the quaternion group
both qualify for $\widetilde \Gamma$ if $\Gamma$ is the Vierergruppe.
\par

To check the existence of stem covers, 
consider $H_2(\Gamma) := H_2(\Gamma, \mathbf Z)$ 
as a trivial $\Gamma$-module,
the short exact sequence
\begin{equation*}
\{0\} 
% \hskip.2cm \longrightarrow \hskip.2cm
\to
\operatorname{Ext}(\Gamma^{\operatorname{ab}}, H_2(\Gamma))
% \hskip.2cm \longrightarrow \hskip.2cm
\to
H^2(\Gamma, H_2(\Gamma))
% \hskip.2cm \longrightarrow \hskip.2cm
\to
\operatorname{Hom}(H_2(\Gamma), H_2(\Gamma))
% \hskip.2cm \longrightarrow \hskip.2cm
\to
\{0\} 
\end{equation*}
of the universal coefficient theorem in cohomology,
and a multiplier $\zeta$ in $Z^2(\Gamma, H_2(\Gamma))$
of which the cohomology class $\underline{\zeta}$
is mapped onto the identity homomorphism
of $H_2(\Gamma)$ to itself.
Then the corresponding central extension
\begin{equation*}
\{0\} 
\hskip.2cm \longrightarrow \hskip.2cm
H_2(\Gamma)
\hskip.2cm \longrightarrow \hskip.2cm
\widetilde \Gamma
\hskip.2cm \overset{p}{\longrightarrow} \hskip.2cm
\Gamma
\hskip.2cm \longrightarrow \hskip.2cm
\{1\} ,
\end{equation*}
in other words the central extension of characteristic class $\underline{\zeta}$, 
is
% \footnote{Wanted: une d\'emonstration directe 
% .......................................................... !!??
% C'est sans doute une cons\'equence de ce qui est rappel\'e
% dans la \emph{footnote} pr\'ec\'edente.} 
a stem cover of $\Gamma$.
Stem covers of $\Gamma$ are classified 
(as central extensions of $\Gamma$ by $H_2(\Gamma)$)
by the group
$\operatorname{Ext}(\Gamma^{\operatorname{ab}}, H_2(\Gamma, \mathbf Z))$;
see Proposition V.5.3 of \cite{Stam--73},
and Theorem 2.2 of \cite{EcHS--72}.
In particular, if $\Gamma$ is perfect, 
it has a unique stem cover, 
also called its \emph{universal central extension}. 
If $\Gamma$ is finite, its stem covers are also called its 
\emph{Schur representation groups}.

Let $p : \widetilde \Gamma \longrightarrow \Gamma$ be a stem cover.
For any central extension
\begin{equation*}
\{0\} 
\hskip.2cm \longrightarrow \hskip.2cm
A
\hskip.2cm \longrightarrow \hskip.2cm
\widetilde \Delta
\hskip.2cm \overset{q}{\longrightarrow} \hskip.2cm
\Delta
\hskip.2cm \longrightarrow \hskip.2cm
\{1\}
\end{equation*}
with divisible kernel $A$ 
(more generally with $A$ such that 
$\operatorname{Ext}(\Gamma^{\operatorname{ab}}, A) = \{0\}$)
and for any homomorphism
$\underline{\rho} : \Gamma \longrightarrow \Delta$,
there exists a homomorphism 
$\rho^0 : \widetilde \Gamma \longrightarrow \widetilde \Delta$
such that $\underline{\rho}(p(\tilde \gamma)) = q(\rho^0 (\tilde \gamma))$
for all $\tilde \gamma \in \widetilde \Gamma$; 
see Proposition V.5.5 of \cite{Stam--73}.  
In particular, for a Hilbert space $\mathcal H$ and a homomorphism 
$\underline{\pi} : \Gamma \longrightarrow \mathcal P \mathcal U (\mathcal H)$,
there exists a unitary representation 
$\pi^0 : \widetilde \Gamma \longrightarrow \mathcal U (\mathcal H)$
such that $\underline{\pi} (p(\tilde \gamma)) = p_{\mathcal H} (\pi^0 (\tilde \gamma))$
for all $\tilde \gamma \in \widetilde \Gamma$.

Observe that, if $\Gamma$ is countable, $H_2(\Gamma)$ is countable
(this follows for example from the Schur-Hopf formula
$H_2(\Gamma) = R \cap [F,F] / [F,R]$ where $\Gamma = F/R$
with $F$ free),
so that $\widetilde \Gamma$ is also countable.

\medskip

It would be interesting to understand, say for the proof of Theorem \ref{Thm4},
if and how one could use the stem cover(s) of $\Gamma$ 
instead of the groups $\Gamma(\zeta)$ 
which appear  in  Section \ref{section4}.

\end{document}